\documentclass[12pt,epsf,bezier]{article}
\renewcommand{\baselinestretch}{1.65}
\usepackage{latexsym}
\usepackage{epsfig}
\usepackage{amssymb}
\usepackage{float}
\usepackage{makeidx}
\usepackage{amsmath}

\begin{document}

\newtheorem{pro}{Proposition}
\newtheorem{theo}{Theorem}
\newtheorem{coro}{Corollary}
\newtheorem{lem}{Lemma}

\newcommand{\bboldy}[1]{\mbox{\boldmath${\#1}$}}
%with this the command $\boldy{\beta}$ equals $\bbeta$ as defined below

\def\E{\mbox{\rm E}}
\def\Var{\mbox{\rm Var}}
\def\Cov{\mbox{\rm Cov}}
\def\Corr{\mbox{\rm Corr}}
\def\Pr{\mbox{\rm Pr}}
\def\I{\mbox{\rm I}}
\def\bbeta{\mbox{\boldmath${\beta}$}}
\def\bgamma{\mbox{\boldmath${\gamma}$}}
\def\sbeta{\mbox{\scriptsize \boldmath${\beta}$}}
\def\eeps{\mbox{\boldmath${\epsilon}$}}
\def\bmu{\mbox{\boldmath${\mu}$}}
\def\bxi{\mbox{\boldmath${\xi}$}}
\def\real{\hbox{\rm\setbox1=\hbox{I}\copy1\kern-.45\wd1 R}}
\newcommand{\bgams}{\mbox{\boldmath{\scriptsize $\gamma$}}}
\newcommand{\bbs}{\mbox{\boldmath{\scriptsize $\beta$}}}
\newcommand{\omt}{\tilde{\omega}}
\newcommand{\bLam}{\mbox{\boldmath{$\Lambda$}}}
\newcommand{\Ell}{\mbox{$\cal{L}$}}
\newcommand{\cG}{\mathcal{G}}
\newcommand{\cL}{\mathcal{L}}
\newcommand{\half}{\frac{1}{2}}
\newcommand{\bsh}{\parindent 0em}
\newcommand{\esh}{\parindent 2.0em}
\newcommand{\eps}{\epsilon}
\newcommand{\peps}{{(\eps)}}
\newcommand{\bzro}{{\bf 0}}
\newcommand{\CP}{\overset{\mathrm{P}}{\rightarrow}}
\newcommand{\cP}{\mathcal{P}}
\def\bU{{\bf U}}
\def\bV{{\bf V}}
\def\bG{{\bf G}}
\def\bC{{\bf C}}
\def\bD{{\bf D}}
\def\bu{{\bf u}}

\begin{center}
\bf \large Case-control survival analysis with a general
semiparametric shared frailty model - a pseudo full likelihood
approach
\end{center}

\normalsize

\vspace*{1cm}

\begin{center}

\bf

Malka Gorfine\footnote{ To whom correspondence should be
addressed} \\ {\it Faculty of Industrial Engineering and
Management, Technion City, Haifa 32000, Israel}\\
{gorfinm@ie.technion.ac.il}
\\

\vspace*{1cm} David M. Zucker\\ {\it Department of Statistics,
Hebrew University, Mt. Scopus, Jerusalem 91905, Israel}\\
{mszucker@mscc.huji.ac.il} \\

\vspace*{1cm} Li Hsu \\ {\it Division of Public Health Sciences,
Fred Hutchinson Cancer Research Center,} \\ {\it Seattle, WA 98109-1024,
USA}\\ {lih@fhcrc.org} \end{center}

\vspace*{1in}

\begin{center}
\today
\end{center}

\newpage

\section*{ Summary}
In this work we deal with correlated failure time (age at onset)
data arising from population-based case-control studies, where
case and control probands are selected by population-based
sampling and an array of risk factor measures is collected for
both cases and controls and their relatives. Parameters of
interest are effects of risk factors on the hazard function of
failure times and within-family dependencies of failure times
after adjusting for the risk factors. Due to the retrospective
nature of sampling, a large sample theory for existing methods has
not been established. We develop a novel estimation techniques for
estimating these parameters under a general semiparametric shared
frailty model. We also present a simple, easily computed, and
non-iterative nonparametric estimator for the cumulative baseline
hazard function. A rigorous large sample theory for the proposed
estimators of these parameters is given along with simulations and
a real data example illustrate the utility of the proposed method.

\section*{ Keywords}
Case-control study; Correlated failure times; Family study;
Frailty model; Multivariate survival model

\newpage

\section{ Introduction}
Clustered failure times arise often in medical and
epidemiologic studies. Examples include disease onset times
of twins (in terms of age),
multiple recurrences times of infections on an individual, or
time to blindness for both eyes within an individual. A typical
case-control family study includes a random sample of independent
diseased individuals (cases) and non-diseased individuals
(controls), along with their family members. An array of genetic
and environmental risk factor measures is collected on these
individuals. Integration of genetic and environmental data is a
central problem of modern observational epidemiology (Hopper et
al., 1994; Zhao et al., 1997; Malone et al., 1998; Malone et al.
2000; Becher et al., 2003). Case-control family studies are
powerful because they provide an efficient way to assess the
effect of risk factors on the occurrence of a rare disease, and
furthermore allow researchers to dissect genetic and environmental
contributions to the disease based on the familial aggregation
pattern of disease clusters. Hopper (2003), in a commentary
article, suggested that such study designs may be the future of
epidemiology, not just genetic epidemiology. Hence the need for
statistical methods that can fully utilize such data is acute.

In this work we focus on population-based case-control family
studies, where a number of case and control probands are randomly
sampled from a well-defined population. The probands are the index
subjects because of whom the families are ascertained. Here we use
the term proband in a broad sense to refer to both cases and
controls, in contrast with the traditional usage in which proband
refers only to cases.

Relative to classical case-control methods, analysis of such
studies is complicated in several ways: (1) Comparisons are no
longer solely between subjects with and without the disease under
study, but rather between collections of the case probands and
their relatives and the control probands and their relatives, each
collection typically including many subjects both with and without
the studied disease. (2) Data are clustered within families, and
hence reflect intra-familial correlation due to unmeasured genetic
and environmental factors.

Our work is motivated by a recent breast cancer study conducted at
the Fred Hutchinson Cancer Research Center (Malone et al., 1998;
Malone et al., 2000). In this study, the cases were incident
breast cancer cases ascertained from the Surveillance,
Epidemiology, and End Results (SEER) registry, which is a set of
geographically defined, population-based cancer registries in the
United States. The controls were selected by random digit dialing,
and were matched with cases based on age at diagnosis and county
of residence. Female relatives of case and control probands were
identified, and the risk factor and outcome information was
subsequently collected on these relatives. The primary goals of
the study are (a) to determine the degree of the strength of the
dependency of ages at diagnosis of breast cancer between probands
and their relatives; (b) assess the effects of covariates on
breast cancer risk.

Two modeling approaches, marginal and conditional, are typically
used for accounting for the correlation within a cluster. In the
conditional model, the correlation is explicitly induced by a
cluster-specific random effect, with the outcomes of the cluster
members being conditionally independent given the random effect.
The random effects model for failure time outcome is generally
known as frailty model, in that the random effect or frailty is
assumed to act multiplicatively on the baseline hazard rate of
failure.  Many frailty models have been considered, including
gamma (Gill, 1985, 1989; Nielsen et al., 1992; Klein 1992, among
others), positive stable (Hougaard, 1986; Fine et al., 2003),
inverse Gaussian, compound Poisson (Aalen, 1992) and log-normal
(McGilchrist, 1993; Ripatti and Palmgren, 2000; Vaida and Xu,
2000, among others). Hougaard (2000) presented a comprehensive
review of the properties of various frailty distributions. Under a
frailty model, the regression coefficients are cluster-specific
log-hazard ratios. By contrast, in the marginal model the
correlation is modelled through a multivariate distribution, such
as a copula function (Genest and MacKay, 1986; Marshall and Olkin,
1988; Shih and Louis, 1995) with a specified model for the
marginal hazard functions. The regression coefficients in the
marginal model represent the log-hazard ratios at the population
level regardless of which cluster an individual comes from. The
effect therefore is ``population-averaged."  Zeger et al. (1988)
provided a comprehensive comparison of the conditional and
marginal modelling approaches.

Methods have been developed for the age at disease onset data from
case-control family studies under both modelling approaches. Shih
and Chatterjee (2002) proposed a semi-parametric
quasi-partial-likelihood approach for estimating the regression
coefficients in a bivariate copula model. Their cumulative hazard
estimator requires an iterative solution, and thus the properties
of their estimators could only be investigated so far by a
simulation study. Moreover, in the presence of multiple relatives
for each proband, the relatives were treated as if they were
independent of each other, which may lead to loss of efficiency in
the baseline hazard function estimator. In contrast, Hsu et al.\
(2004) presented a quasi-EM algorithm method for the popular gamma
frailty model. In the random effects model of Hsu et al., the
regression coefficients express the effect on a subject's disease
risk due to being exposed relative to the same subject's level of
risk when unexposed.  The baseline hazard function estimator
naturally accommodates multiple relatives in a family (Hsu and
Gorfine, 2006).  However, the properties of the proposed
estimators were also studied only by simulation. The method of
Shih and Chatterjee (2002) can be adapted to the family-specific
frailty setting (Oakes, 1989), but with the same limitation as for
the marginal model: the lack of large sample theory.
% and a potentially
%computationally slow iterative procedure for the baseline hazard
%function estimators.

In this work, we develop a new estimation technique for the general
semiparametric shared frailty model, where the parameters of
interest are the regression coefficients and the frailty
parameters. Our general family-specific frailty model is for any
frailty distribution that has finite moments.
%The proposed
%baseline hazard function estimator is non-iterative which is easy
%in implementation and fast in computation.
The estimation procedure for the baseline hazard function leads to
an estimator whose asymptotic properties can be derived and
expressed in a tractable manner.
%The details of the asymptotic theory, involving a proof that the
%estimators of the regression coefficients and the frailty
%parameters are consistent and asymptotically normal, along with a
%consistent closed-form estimator for the asymptotic covariance
%matrix of the parameter estimates, are beyond the scope of this
%paper and will be presented in a separate communication.

Section 2 presents our model, and Section 3 describes our estimation
procedure. Section 4 gives the consistency and asymptotic normality
results for the estimators. In Section 5, we describe an extension
of our method for the case where the proband observation times are
subject to a certain restriction that can arise in some studies.
Section 6 presents simulation results. In Section 7 we illustrate our
method with a case-control family study of early onset breast cancer.
Section 8 provides a short discussion. The Appendix provides the
details of the asymptotic theory.

\section{ Notation and model formulation}
We consider a matched case-control family study where one case
proband is age-matched with one control proband, and an array of
risk factors is measured on the case and control probands and
their relatives. Each matched set contains one case family and one
control family, and there are $n$ i.i.d.\ matched sets. Let
$T^0_{ij}$ and $C_{ij}$ denote the age of disease onset and age at
censoring, respectively, for individual $j$ of family $i$,
$i=1,\ldots,2n$, $j=0,1,\ldots,m_i$, where $j=0$ corresponds to
the proband. Following Parner (1998, p.~187), we regard $m_i$ as a
random variable over $\{1, \ldots, m \}$ for some $m$, and build
up the remainder of the model conditional on $m_i$. Define
$\delta_{ij}=I(T^0_{ij} \leq C_{ij})$ to be the failure indicator
and $T_{ij}=\min(T^0_{ij},C_{ij})$ to be the observed follow-up
time for individual $ij$. We assume that a $p$-vector of
covariates is observed on all subjects, and let ${\bf Z}_{ij}$
denote the value of the (time-independent) covariate vector for
individual $ij$. In addition, we associate with family $i$ an
unobservable family-level covariate $\omega_i$, the ``frailty",
which induces dependence among family members. The conditional
hazard function for proband $i$, given the family frailty
$\omega_i$, is assumed to take the form
 \begin{equation}\label{eq:haz2}
\lambda_{i0}(t|{\bf Z}_{i0},\omega_i)=\omega_i \lambda_0(t)
\exp(\bbeta^T {\bf Z}_{i0}) \;\;\;\;\;\; i=1,\ldots,2n.
 \end{equation}
The conditional hazard function for relative $ij$,
$j=1,\ldots,m_i$, given the family frailty $\omega_i$ and proband
$i$'s data, is assumed to take the form
 \begin{equation}\label{eq:haz1}
\lambda_{ij}(t|T_{i0},\delta_{i0},{\bf Z}_{i0},{\bf Z}_{ij},
\omega_i)=\omega_i \lambda_0(t) \exp(\bbeta^T {\bf Z}_{ij})
\;\;\;\;\;\; i=1,\ldots,2n; \;\;\; j=1,\ldots,m_i.
 \end{equation}
Here $\bbeta$ is a $p$-vector of unknown regression coefficients,
and $\lambda_0$ is a conditional baseline hazard of unspecified
form. The above model implies that the proband and the relatives have a
common conditional baseline hazard function $\lambda_0$, and that
all the dependence between the proband and the relatives in a
given family is due to the frailty factor $\omega_i$. The random variable
$\omega_i$ is assumed to have a density $f(\omega)\equiv f(w;\theta)$,
where $\theta$ is an unknown parameter.
For simplicity, we assume that $\theta$ is a scalar, though the vector
case could be developed in a similar manner.

We put $\bgamma=
({\bbeta}^T,\theta)^T$, and let $\bgamma^\circ=({\bbeta^\circ}^T,
\theta^\circ)^T$ denote the true value of $\bgamma$.
The objective is to
estimate $\bgamma$ and $\Lambda_0(t)=\int_0^t
\lambda_0(u)du$. Let $\Lambda_0^\circ(t)$ denote the true
value of $\Lambda_0$. Further, let $\delta_{iR}=(\delta_{i1},\ldots,\delta_{im_i})$,
$T_{iR}=(T_{i1},\ldots,T_{im_i}),$ and
${\bf Z}_{iR}=({\bf Z}_{i1},\ldots,{\bf Z}_{im_i})$.

We make the following basic assumptions.
\begin{enumerate}
\item ${\bf Z}_{ij}$ is bounded. \item The parameter $\bgamma$
lies in a compact subset $\cG$ of $\real^{p+1}$ containing an open
neighborhood of $\bgamma^\circ$. \item Conditional on $\{{\bf
Z}_{ij}\}_{j=1}^{m_i}$ and $\omega_i$, the censoring times are
independent and noninformative for $\omega_i$ and $(\bbeta,
\Lambda_0)$. In addition, the frailty $\omega_i$ is independent of
$\{{\bf Z}_{ij}\}_{j=1}^{m_i}$. \item The effect of the covariates
on age at onset is subject-specific, i.e.\
$\Pr(T_{ij},\delta_{ij}|{\bf Z}_{i0},{\bf Z}_{iR},\omega_i)$
$=\Pr(T_{ij},\delta_{ij}|{\bf Z}_{ij}, \omega_i)$. This implies
$\Pr(T_{ij},\delta_{ij}|{\bf Z}_{i0},{\bf
Z}_{iR})=\Pr(T_{ij},\delta_{ij}|{\bf Z}_{ij})$.
\end{enumerate}

The first two of these assumptions imply that there exists a positive
constant $\nu$ such that
\begin{equation}
\nu^{-1} \leq \exp(\bbeta^T {\bf Z}_{ij}) \leq \nu.
\label{eb}
\end{equation}
A number of additional technical assumptions are listed in the appendix.

The likelihood function for the data can be written as
\begin{eqnarray}\label{eq:like}
\lefteqn{L}&=&\prod_{i=1}^{2n}
 f(T_{iR},\delta_{iR},{\bf Z}_{iR},{\bf Z}_{i0}|T_{i0},\delta_{i0}) \nonumber \\
 &=& \prod_{i=1}^{2n}
 f(T_{iR},\delta_{iR}|{\bf Z}_{iR},{\bf Z}_{i0},T_{i0},\delta_{i0})
\times f({\bf Z}_{iR}|{\bf Z}_{i0}) \times
 f({\bf Z}_{i0}|T_{i0},\delta_{i0}).
\end{eqnarray}
Since $f({\bf Z}_{iR}|{\bf Z}_{i0})$ does not depend on the
parameters of interest $(\bbeta,\Lambda_0,\theta)$, this
term will be ignored. In the following subsections we consider the
other two terms in (\ref{eq:like}).

\subsection{ The likelihood for the proband data}
For the likelihood function of the proband data, $\prod_i^{2n}
f({\bf Z}_{i0}|T_{i0},\delta_{i0})$, we use a retrospective
likelihood for the standard case-control study (Prentice and
Breslow, 1978). We express this likelihood in terms of the
marginal survival function $
 {S_{i0}(t)}= Pr(T_{i0}>t|{\bf Z}_{i0})
 =\int
Pr(T_{i0}>t|{\bf Z}_{i0},\omega)f(\omega) d\omega.
 $
 In our
setting we have $n$ one-to-one matched sets. Based on the marginal
survivor function, the marginal hazard function can be written as
 $$
\lambda_{i0}(t)=\lambda_0(t) \exp(\bbeta^T {\bf Z}_{i0})
\frac{\mu_{1i}(t;\bgamma,\Lambda_0)}{\mu_{0i}(t;\bgamma,\Lambda_0)},
 $$
where
 $$
\mu_{ki}(t;\bgamma,\Lambda_0) =
 \int \omega^k \exp\{- \omega H_{i0}(t) \} f(\omega) d \omega
 \;\;\;\;\;\;\;\;\;\;  k=0,1,2
 $$
and
 $
H_{i0}(t) = \Lambda_0(t)\exp(\bbeta^T {\bf Z}_{i0}).
 $
We arrange the notation so that the first $n$ families are the
case families and the $r$th case family, $r=1,\ldots,n$, is
matched with the $(n+r)$th control family. The likelihood for
the proband data is then replaced by the following conditional
likelihood:
\begin{equation}\label{eq:problik}
L^{(1)}=\prod_{r=1}^n \frac{\exp(\bbeta^T {\bf Z}_{r0})
{\xi_{10r}(T_{r0};\bgamma,\Lambda_0)}}
 {\exp(\bbeta^T {\bf Z}_{r0})
 {\xi_{10r}(T_{r0};\bgamma,\Lambda_0)}
+
 \exp(\bbeta^T {\bf Z}_{(n+r)0})
 {\xi_{10(n+r)}(T_{(n+r)0};\bgamma,\Lambda_0)}},
\end{equation}
where
 $$
 \xi_{kk'i}(t;\bgamma,\Lambda_0)=
 \frac{\mu_{ki}(t;\bgamma,\Lambda_0)}{\mu_{k'i}(t;\bgamma,\Lambda_0)},
\;\;\;\;\;\;\;\;\;\;\;\;\; k,k'=0,1,2.
 $$
Let
 $$
 \xi_{10i}^{\beta_l}(t;\bgamma,\Lambda_0)= \frac{\partial}{\partial \beta_l}
 \xi_{10i}(t;\bgamma,\Lambda_0)
=H_{i0}(t)Z_{i0l}\left\{\xi_{10i}^2(t;\bgamma,\Lambda_0)-\xi_{20i}(t;\bgamma,\Lambda_0)
\right\}
 ,
 $$
 $$
 \mu_{ki}^{\theta}(t;\bgamma,\Lambda_0)=\frac{\partial}{\partial \theta}
 \mu_{ki}(t;\bgamma,\Lambda_0)=
 \int \omega^k \exp\{-\omega H_{i0}(t)\}\frac{\partial}{\partial
 \theta}f(\omega)d\omega,
 $$
 and
 \begin{eqnarray}
\xi_{kk'i}^{\theta}(t;\bgamma,\Lambda_0)&=&
\frac{\partial}{\partial \theta}
 \xi_{kk'i}(t;\bgamma,\Lambda_0)  \nonumber \\
&=&
\frac{\mu_{k'i}(t;\bgamma,\Lambda_0)\mu_{ki}^{\theta}(t;\bgamma,\Lambda_0)
-\mu_{k'i}^{\theta}(t;\bgamma,\Lambda_0)\mu_{ki}(t;\bgamma,\Lambda_0)}
{\mu_{k'i}^2(t;\bgamma,\Lambda_0)}, \nonumber
 \end{eqnarray}
 for $k,k'=0,1,2$ and $l=1,\ldots,p$. Then the score function
 for $\beta_l$, $l=1,\ldots,p$,
 is given by
\begin{eqnarray}
U^{(1)}_l(\bgamma,\Lambda_0)&=&
 \sum_{r=1}^n \left\{ {\bf Z}_{r0}+\frac{\xi_{10r}^{\beta_l}(T_{r0};\bgamma,\Lambda_0)}
 {\xi_{10r}(T_{r0};\bgamma,\Lambda_0)} \right. \nonumber \\
&-& \frac{\exp(\bbeta^T {\bf
Z}_{r0})[Z_{r0l}\xi_{10r}(T_{r0};\bgamma,\Lambda_0)+\xi_{10r}^{\beta_l}
 (T_{r0};\bgamma,\Lambda_0)]}
 {\exp(\bbeta^T {\bf Z}_{r0})\xi_{10r}(T_{r0};\bgamma,\Lambda_0)
 +\exp(\bbeta^T {\bf
 Z}_{(n+r)0})\xi_{10(n+r)}(T_{(n+r)0};\bgamma,\Lambda_0)} \nonumber\\
&-& \left. \frac{\exp(\bbeta^T {\bf
Z}_{(n+r)0})[Z_{(n+r)0l}\xi_{10(n+r)}(T_{(n+r)0};\bgamma,\Lambda_0)
 +\xi_{10(n+r)}^{\beta_l}
 (T_{(n+r)0};\bgamma,\Lambda_0)]}
 {\exp(\bbeta^T {\bf Z}_{r0})\xi_{10r}(T_{r0};\bgamma,\Lambda_0)
 +\exp(\bbeta^T {\bf
 Z}_{(n+r)0})\xi_{10(n+r)}(T_{(n+r)0};\bgamma,\Lambda_0)}
 \right\}, \nonumber
\end{eqnarray}
and the score function for $\theta$ is given by
\begin{eqnarray}
U^{(1)}_{p+1}(\bgamma,\Lambda_0) &=& \sum _{r=1}^n \left\{
\frac{\xi_{10r}^{\theta}(T_{r0};\bgamma,\Lambda_0)}{\xi_{10r}(T_{r0};\bgamma,\Lambda_0)}
\right.\nonumber\\
 &-&\left.
 \frac{\exp(\bbeta^T {\bf Z}_{r0})\xi_{10r}^{\theta}(T_{r0};\bgamma,\Lambda_0)+
 \exp(\bbeta^T
  {\bf Z}_{(n+r)0})\xi_{10(n+r)}^{\theta}(T_{(n+r)0};\bgamma,\Lambda_0)}
  {\exp(\bbeta^T {\bf Z}_{r0})\xi_{10r}(T_{r0};\bgamma,\Lambda_0)+
 \exp(\bbeta^T
  {\bf
  Z}_{(n+r)0})\xi_{10(n+r)}(T_{(n+r)0};\bgamma,\Lambda_0)}\right\}.
  \nonumber
\end{eqnarray}
Under the gamma frailty model, we have
\begin{eqnarray*}
\frac{\mu_{1i}(t;\bgamma,\Lambda_0)}{\mu_{0i}(t;\bgamma,\Lambda_0)}=
\left\{\theta H_{i0}(t)+1 \right\}^{-1},
\end{eqnarray*}
and so the likelihood function (\ref{eq:problik}) corresponds to that
presented in Hsu et al.\ (2004) in the case of one-to-one
matching. Extension to matching of multiple cases or multiple
controls are straightforward, see e.g. Breslow and Day (1980).

\subsection{ The likelihood for the data from the relatives}
Let $N_{ij}(t)=\delta_{ij}I(T_{ij}\leq t)$, $j=1,\ldots,m_i$,
$N_{i.}(t)=\sum_{j=1}^{m_i}N_{ij}(t)$, $H_{ij}(t)=\Lambda_0(T_{ij}
\wedge~t)$ $\exp(\bbeta^T {\bf Z}_{ij})$, $j=1,\ldots,m_i$, and
$H_{i.}(t)=\sum_{j=1}^{m_i}H_{ij}(t)$, and let $\tau$ be the
maximum follow-up time. The likelihood of the data from the
relatives then can be written as
\begin{eqnarray}
L^{(2)}&=&\prod_{i=1}^{2n} \int \prod_{j=1}^{m_i}
 \left\{ \lambda_{ij}(T_{ij}|T_{i0},\delta_{i0},{\bf
 Z}_{i0},\omega)\right\}^{\delta_{ij}} S_{ij}(T_{ij}|T_{i0},\delta_{i0},{\bf
 Z}_{i0},\omega)f(w|T_{i0},\delta_{i0},{\bf
 Z}_{i0})dw \nonumber \\
 &=& \prod_{i=1}^{2n} \prod_{j=1}^{m_i}
 \left\{\lambda_0(T_{ij})\exp(\bbeta^T {\bf Z}_{ij}) \right\}^{\delta_{ij}}
 \prod_{i=1}^{2n} \int \omega^{N_{i.}(\tau)} \exp\{-\omega H_{i.}(\tau)\}
 f(\omega|T_{i0},\delta_{i0},{\bf Z}_{i0})d\omega. \nonumber
\end{eqnarray}
Here, by a Bayes theorem argument,
\begin{equation}
f(\omega|T_{i0},\delta_{i0},{\bf Z}_{i0})= \frac
{\omega^{\delta_{i0}} \exp(-\omega \Lambda_0(T_{i0}) e^{\bbs^T {\bf Z}_{i0}}) f(\omega)}
{\int \omt^{\delta_{i0}} \exp(-\omt \Lambda_0(T_{i0}) e^{\bbs^T {\bf Z}_{i0}}) f(\omt) \, d\omt} \, .
\label{cdo}
\end{equation}
The log-likelihood is given by
 $$
l^{(2)}= \sum_{i=1}^{2n} \sum_{j=1}^{m_i} \delta_{ij}
\log\{\lambda_0(T_{ij})\exp(\bbeta^T {\bf Z}_{ij}) \}
 + \sum_{i=1}^{2n}\log \left\{\int \omega^{N_{i.}(\tau)} \exp\{-\omega H_{i.}(\tau)\}
 f(\omega|T_{i0},\delta_{i0},{\bf Z}_{i0})d\omega \right\}.
 $$
The scores for $(\beta_1,\ldots,\beta_p)$ are given by
 \begin{eqnarray}
U^{(2)}_l(\bgamma,\Lambda_0)&=&\sum_{i=1}^{2n} \sum_{j=1}^{m_i}
\delta_{ij} Z_{ijl}
 +  \sum_{i=1}^{2n} \frac{\int \omega^{N_{i.}(\tau)} \exp\{-\omega
H_{i.}(\tau)\}
 \frac{\partial}{\partial \beta_l} f(\omega|T_{i0},\delta_{i0},{\bf Z}_{i0})d\omega }
 {\int \omega^{N_{i.}(\tau)} \exp\{-\omega H_{i.}(\tau)\}
 f(\omega|T_{i0},\delta_{i0},{\bf Z}_{i0})d\omega } \nonumber \\
 &-&
\sum_{i=1}^{2n} \frac{\int \omega^{N_{i.}(\tau)+1} \exp\{-\omega
H_{i.}(\tau)\}
  f(\omega|T_{i0},\delta_{i0},{\bf Z}_{i0})d\omega }
 {\int \omega^{N_{i.}(\tau)} \exp\{-\omega H_{i.}(\tau)\}
 f(\omega|T_{i0},\delta_{i0},{\bf Z}_{i0})d\omega }
 \sum_{j=1}^{m_i}H_{ij}(\tau)Z_{ijl} \nonumber
 \end{eqnarray}
for $l=1,\ldots,p$. The score for $\theta$ is given by
 $$
 U^{(2)}_{p+1}(\bgamma,\Lambda_0)=\sum_{i=1}^{2n} \frac{\int \omega^{N_{i.}(\tau)} \exp\{-\omega
H_{i.}(\tau)\}
\frac{\partial}{\partial \theta} f(\omega|T_{i0},\delta_{i0},{\bf Z}_{i0})d\omega}{\int \omega^{N_{i.}(\tau)}
\exp\{-\omega H_{i.}(\tau)\}
f(\omega|T_{i0},\delta_{i0},{\bf Z}_{i0})d\omega}.
$$

\section { The proposed approach}
We focus first on estimating the baseline cumulative hazard
function $\Lambda_0(t)$. Let $Y_{ij}(t)=I(T_{ij} \geq t)$, and let
${\mathcal F}_t$ denote the $\sigma$-algebra generated by
$(T_{i0},\delta_{i0},{\bf Z}_{i0})$ plus the entire observed
history of the relatives up to time $t$:
 $$
{\mathcal F}_{t} = \sigma(T_{i0},\delta_{i0},{\bf Z}_{i0},
N_{ij}(u), Y_{ij}(u), {\bf Z}_{ij}; i=1,\ldots,2n;
j=1,\ldots,m_i;0 \leq u \leq t).
 $$
Then, as discussed by Gill (1992) and Parner (1998),  the
stochastic intensity process for $N_{ij}(t)$, $i=1,\ldots,2n$,
$j=1,\ldots,m_i$, with respect to ${\mathcal F}_{t}$ is given by
\begin{eqnarray}\label{eq:inten}
\lambda_0(t)\exp(\bbeta^T {\bf Z}_{ij}) Y_{ij}(t)
\psi_{i}(t-,\bgamma,\Lambda_0),
\end{eqnarray}
where, using (\ref{cdo}),
\begin{eqnarray}
\psi_{i}(t,\bgamma,\Lambda_0) &=& \E[\omega_i|{\mathcal
F}_{t}] \nonumber \\
 &=& \frac{\int \omega^{N_{i.}(t)+1}\exp(-\omega H_{i.}(t))f(\omega|T_{i0},\delta_{i0},{\bf Z}_{i0})d\omega}
 {\int \omega^{N_{i.}(t)}\exp(-\omega H_{i.}(t))f(\omega|T_{i0},\delta_{i0},{\bf
 Z}_{i0})d\omega} \nonumber \\
 &=& \frac{\int \omega^{N_{i.}(t)+1+\delta_{i0}}\exp(-\omega \{H_{i.}(t) + H_{i0}(T_{i0})\}) f(\omega) d\omega}
{\int \omega^{N_{i.}(t)+\delta_{i0}}\exp(-\omega \{H_{i.}(t) + H_{i0}(T_{i0}) \}) f(\omega) d\omega} \nonumber
\, .
\end{eqnarray}

Define
(for $0 \leq r \leq m$ and $h \geq 0$)
\begin{equation}
\psi^*(r,h) = \frac{\int
w^{r+1} e^{-hw} f(w) dw} {\int w^{r} e^{-hw} f(w) dw}.
\label{psi}
\end{equation}
Some salient properties of $\psi^*(r,h)$ are noted in Sec.~9.2.
With this definition, we have $\psi_{i}(t,\bgamma,\Lambda_0) =
\psi^*(N_i(t),H_{i \cdot}(t))$.

The key to obtaining parameter estimators for a semiparametric
survival model is an estimator of the nonparametric baseline
hazard function. For our model, a Breslow-type estimator with a
jump at each observed failure time among the relatives can be
formulated in a natural way (Shih and Chatterjee, 2002). However,
the hazard function for the relatives at time $t$ depends on their
respective proband's observation time $T_{i0}$. For example, under
the gamma frailty model with expectation 1 and variance $\theta$,
$\psi_i(t,\bgamma,\Lambda_0)=\{\theta^{-1}+N_{i.}(t)+\delta_{i0}
\}\{\theta^{-1}+H_{i.}(t)+H_{i0}(T_{i0})\}^{-1}$. Often the
relevant proband's observation time is greater than $t$, so that
the standard Breslow formula for the baseline hazard estimator at
time $t$ involves values of $\Lambda_0$ for times beyond time $t$.
An iterative procedure is thus required to obtain the estimator.
In addition, because of this estimator's complicated structure,
its asymptotic properties have not been established.

We propose to estimate the baseline hazard function using a non-iterative
two-stage procedure. The first-stage estimator is a weighted Breslow-type
estimator, where the weight at time $t$ for family $i$ is equal to 1 if the
observation time $T_{i0}$ of the family $i$ proband is less that $t$, and
equal to 0 otherwise.
The
second-stage estimator is the standard Breslow-type estimator that
uses all the relatives' failure times, plugging in the first-stage
estimator where necessary.

More specifically, the estimators are defined as follows. Let $\tau_g$,
$g=1,\ldots,G$, denote the observed failure times of the relatives
and assume that $d_g$ failures were observed at time $\tau_g$.
In theory, since we are dealing with continuous survival distributions,
$d_g=1$ for all $g$, but we write the formula for the estimator in a
form that allows for a modest level of ties in the survival times.
Let $\Lambda_{max}$ be some known (possibly large) upper bound for
$\Lambda_0^\circ(t)$. Define $\bar{\psi}(r,h)= \psi^*(r,h \wedge h_{max})$,
with $h_{max}$ = $m \nu \Lambda_{max}$, where $\nu$ is as in (\ref{eb}).
Further, define $\bar{\psi}_i(t,\bgamma,\Lambda) = \bar{\psi}
(N_i(t),H_{i \cdot}(t,\bgamma,\Lambda))$.
The first-stage estimator is then defined as a step function
whose $g$-th jump is given by
\begin{equation}
\label{eq:stage1}
\Delta \tilde{\Lambda}_0(\tau_g)=
 \frac {\sum_{i=1}^{2n}I(T_{i0}<\tau_g)
\sum_{j=1}^{m_i}dN_{ij}(\tau_g)} {\sum_{i=1}^{2n}I(T_{i0}<\tau_g)
\bar{\psi}_i(\tau_{g-1},\bgamma,\tilde{\Lambda}_0)
 \sum_{j=1}^{m_i}Y_{ij}(\tau_g)\exp(\bbeta^T {\bf Z}_{ij})}.
 \end{equation}
In a similar way, the second-stage estimator is defined as a step function
whose $g$-th jump is given by
 \begin{equation}
 \label{eq:stage2}
 \Delta
\hat{\Lambda}_0(\tau_g)=
 \frac {d_g}
{\sum_{i=1}^{2n}\tilde{\psi}_i(\tau_{g-1},\bgamma)
 \sum_{j=1}^{m_i}Y_{ij}(\tau_g)\exp(\bbeta^T {\bf Z}_{ij})},
 \end{equation}
where $\tilde{\psi}_i(t,\bgamma)$ is defined
analogously to $\bar{\psi}_i(t,\bgamma,{\Lambda}_0)$, with
$\Lambda_0(T_{i0})$ replaced by $\tilde\Lambda_0(T_{i0})$ if
$T_{i0}\geq t$ and by $\hat{\Lambda}_0(T_{i0})$ otherwise. It is
clear that no iterative optimization process is required here and
the large-sample properties of $\hat{\Lambda}_0(t)$ will be
determined by those of $\tilde{\Lambda}_0(t)$.

We note that there is no guarantee that $\tilde{\Lambda}_0(t,\bgamma)$ as defined
above will be bounded by $\Lambda_{max}$, but this does not matter: if desired, we can
replace the estimator by $\min\{\tilde{\Lambda}_0(t,\bgamma), \Lambda_{max}\}$ without
affecting the asymptotics.

For estimating $(\bbeta,\theta)$ we use a pseudo-likelihood
approach: in the score functions based on $L^{(1)}$ and $L^{(2)}$,
we replace the unknown $\Lambda_0$ by $\hat{\Lambda}_0$. Thus, the
score function corresponding to $\beta_l$ (for $l=1,\ldots,p$) is given by $
 U_l(\bgamma,\hat{\Lambda}_0)=n^{-1} \left\{U^{(1)}_l(\bgamma,\hat{\Lambda}_0)
 +U^{(2)}_l(\bgamma,\hat{\Lambda}_0)\right\},
$ and the estimating function for $\theta$ is given by
 $
U_{p+1}(\bgamma,\hat{\Lambda}_0)=n^{-1}
\left\{U^{(1)}_{p+1}(\bgamma,\hat{\Lambda}_0)
+U^{(2)}_{p+1}(\bgamma,\hat{\Lambda}_0)\right\}.
 $
To summarize, our proposed estimation procedure is as follows:
\begin{enumerate}
 \item Provide an initial value for $\bgamma$.
 \item For the given values of $\bgamma$, estimate $\Lambda_0$
using (9) and (10).
 \item For the given value of
$\Lambda_0$, estimate $\bgamma$.
 \item Repeat Steps 2
and 3 until convergence is reached with respect to
$\hat{\Lambda}_0$ and $\hat{\bgamma}$.
\end{enumerate}

\section{ Asymptotic properties}

We show that $\hat{\bgamma}$
is a consistent estimator of $\bgamma^\circ$ and that $\sqrt{n} (\hat{\bgamma}-\bgamma^\circ)$ is
asymptotically mean-zero multivariate normal. In this section, we present a broad
outline sketch of the argument. The Appendix provides the details of the proofs, including
a detailed list of the technical conditions assumed. The arguments are patterned after those
of Gorfine et al.\ (2006) and Zucker et al.\ (2006), but with considerable expansion.

Consistency is shown through the following steps.
\begin{description}
 \item[Claim A1.] $\tilde{\Lambda}_0(t,\bgamma)$ converges in probability to
some function $\Lambda_0^*(t,\bgamma)$ uniformly in $t$ and $\bgamma$.
The
function $\Lambda_0^*(t,\bgamma)$ satisfies $\Lambda_0^*(t,\bgamma^\circ)=\Lambda_0^\circ(t)$.
 \item[Claim A2.] $\hat{\Lambda}_0(t,\bgamma)$ converges in probability to
some function $\Lambda_0(t,\bgamma)$ uniformly in $t$ and $\bgamma$.
The
function $\Lambda_0(t,\bgamma)$ satisfies $\Lambda_0(t,\bgamma^\circ)=\Lambda_0^\circ(t)$.
 \item[Claim B.]
${\bf U}(\bgamma,\hat{\Lambda}_0(\cdot,\bgamma))$ converges in probability
uniformly in $t$ and $\bgamma$ to a limit ${\bf u}
(\bgamma,\Lambda_0(\cdot,\bgamma))$.
 \item[Claim C.] There exists a
unique consistent (in pr.) root to ${\bf
U}(\hat{\bgamma},\hat{\Lambda}_0(\cdot,\hat{\bgamma}))=\bzro$.
\end{description}
The proofs of Claims A1, A2, and B involve empirical process and function-space compactness
arguments, while Claim C is shown using Foutz's (1977) theorem on consistency of maximum
likelihood type estimators.

The proof of asymptotic normality is based on the following decomposition:
\begin{eqnarray}
\lefteqn{\bzro} &=& {\bf U}
(\hat{\bgamma},\hat{\Lambda}_0(\cdot,\hat{\bgamma}))\nonumber
\\ &=& {\bf U}(\bgamma^\circ,\Lambda_0^\circ) +
[{\bf U} (\bgamma^\circ,\hat{\Lambda}_0(\cdot,\bgamma^\circ))-{\bf
U} (\bgamma^\circ,\Lambda_0^\circ)] \nonumber\\ & & + \, [{\bf
U}(\hat{\bgamma},\hat{\Lambda}_0(\cdot,\hat{\bgamma}))- {\bf
U}(\bgamma^\circ,\hat{\Lambda}_0(\cdot,\bgamma^\circ))]. \nonumber
\end{eqnarray}
In the Appendix we analyze each of the above three terms and prove that
$\sqrt{n}(\hat{\bgamma}-\bgamma^\circ)$ has an asymptotic mean-zero multivariate
normal distribution. Although it is possible to develop a consistent
closed-form sandwich estimator for the asymptotic covariance matrix of this
distribution, we do not present this estimator because it is too complicated
to be practically useful. Instead, as discussed in Section 6, we recommend
bootstrap standard-error estimates.

\section{ Extension to restricted sampling of probands}

A key assumption in our procedure for estimating ${\Lambda}_0$ is that the support
of the proband observation times and that of relatives' observation times have the
same lower limit, which is designated (without loss of generality) as time zero.
In some applications, however,
the probands' observed times are restricted to some range $[s_0, s_1]$ with $s_0>0$.
For example, Malone et al.\
(2006) present a multi-center case-control breast cancer study
where ages of cases and controls are restricted between ages
35-64. In a design of this form, where the probands' observed
times are left-restricted by $s_0$ and the relatives' failure
times are unrestricted, $\Lambda_0$ will be underestimated by
our two-stage procedure. But this bias can be easily corrected by
first estimating $\Lambda_0(s_0)$.

We present here the resulting three-stage estimator for the
left-restricted design. Let
$\Delta\tilde{\Lambda}_0\{\tau_g,\Lambda_0(s_0)\}$ and
$\Delta\hat{\Lambda}_0\{\tau_g,\Lambda_0(s_0)\}$ be defined
analogously to $\Delta\tilde{\Lambda}_0(\tau_g)$ and
$\Delta\hat{\Lambda}_0(\tau_g)$ with
$\Lambda_0(T_{i0})=\Lambda_0(s_0)+\sum_{\tau_g \in
[s_0,T_{i0}]}\Delta_0(\tau_g)$.  The estimator
$\hat\Lambda_0(s_0)$ is defined to be the root of
\begin{equation}
\sum_{\tau_g \in
[0,s_0]}\Delta\hat{\Lambda}_0\{\tau_g,\Lambda_0(s_0)\}-\Lambda_0(s_0)=0.
\end{equation}
The root can be found by simple univariate Newton-Raphson iteration. This
completes the first stage. The second stage involves calculating $\Delta\tilde\Lambda_0
\{\tau_g,\hat\Lambda_0(s_0)\}$, $g=1,\ldots, G$, using the formula (\ref{eq:stage1}).
In the third stage, we use the results of the second stage and the formula (\ref{eq:stage2})
to calculate the the final estimate $\Delta\hat\Lambda_0(\tau_g)$, $g=1,\ldots,G$. In applying
(\ref{eq:stage2}),
we replace $\Lambda_0(T_{i0})$ by $\tilde\Lambda_0\{T_{i0},\hat\Lambda_0(s_0)\}$
if $T_{i0}\geq \tau_g$ and by $\hat{\Lambda}_0(T_{i0})$ otherwise.

In Section 6 below, we present simulation results for this estimator. In theory, the
asymptotic properties of the three-stage procedure could be worked out via an extension
of the arguments for the two-stage procedure, but the algebra becomes very complicated.
We hope to develop asymptotic theory for the left-restricted design in future work.

\section { Simulation results - gamma frailty}

We have performed a simulation study to evaluate the finite sample
performance of the proposed method and compare it with existing
methods. One of the most extensively used frailty models is the
model with gamma-distributed frailty. Under this model, $\theta$
quantifies the heterogeneity of risk among families. The larger
the value of $\theta$ is, the stronger the dependence among family
members. In addition, the gamma frailty model can be re-expressed
in terms of the Clayton-Oakes copula-type model (Clayton, 1978;
Oakes, 1989) and the cross-ratio, introduced by Oakes (1989) as a
local measure of association between survival times, is constant
on the support of failure time region and equals $1+\theta$. The
gamma frailty model is also convenient mathematically, because it
admits a closed-form representation of the marginal survival
distributions. These features make the gamma frailty model very
popular. Hence we conducted our simulation study under the gamma
frailty model, using, as is customary, the gamma distribution with
expectation 1 and variance $\theta$.

Simulation results are based on 500 control probands matched to
500 case probands, with one relative sampled for each proband. We
considered a single $U[0,1]$ distributed covariate with
$\beta=\ln(2)$, $\Lambda_0(t)=t$, $\theta=2$, and a $U[0,1]$
censoring  variable, yielding a censoring rate among the relatives
of approximately $60\%$. In Table 1 we compare the following three
estimates: the proposed estimate with the two-stage procedure for
$\Lambda_0$, the estimate of Hsu et al. (2004), and a modified
version of Shih and Chatterjee's (2002) estimate, with their
method adapted to the gamma frailty model. Results are based on
500 simulated data sets. The efficiency difference between our
two-stage estimator and that of Shih and Chatterjee is very small.

For our estimators, in addition to the above-mentioned simulation setting,
we also considered $\beta=0$, $\theta=3$ and a censoring distribution of $U[0,4]$
with a $U[0,4]$ distributed covariate, or a censoring distribution
of $U[0,0.1]$ with $U[0,1]$ distributed covariate, yielding
censoring rates of approximately 30\% or 90\%, respectively.  To
construct confidence intervals, we use a bootstrap approach. In
the setting of censored survival data, the usual nonparametric
bootstrap is problematic because it leads to a substantial
proportion of tied survival times. Hence we used the weighted
bootstrap approach of Kosorok et al. (2004) instead. For the
weighted bootstrap, a sample of $2n$ independent and identically
distributed weights from the unit exponential distribution  was
generated for each bootstrap sample. Let $\xi_1,\ldots,\xi_{2n}$
be the standardized weights after dividing each weight by the
average weight. Then, in the estimating functions, for any given
function $h$ the empirical mean $n^{-1}\sum_{i=1}^{2n}
h(T_i,\delta_i,{\bf Z }_i)$ is replaced by its corresponding
weighted empirical mean $n^{-1}\sum_{i=1}^{2n} \xi_i
h(T_i,\delta_i,{\bf Z }_i)$. Kosorok et al. (2004) proved that
this weighted bootstrap procedure gives valid inference for all
parameters under right-censored univariate failure times.

Results based on the two-stage procedure for $\Lambda_0$ are
presented in Tables 2-4 for various levels of censoring.
We present the mean, the empirical standard
error, and the coverage rate of the $95\%$ weighted bootstrap
confidence interval. The results are based on 50 bootstrap samples
for each of the 2000 simulated data sets of each configuration.
Our estimates perform well in terms of bias and coverage
probability.

For studying the case of left-restricted data, we considered a
similar configuration as of Table 1, but now the probands
observation times are restricted to be $>0.1$. In Table 5 the
results of our three-stage estimator are presented along with the
estimators of Hsu et al. (2004) and of Shih and Chatterjee (2002).
It is seen that estimating $\Lambda_0(s_0)$ yields small
efficiency loss in $\hat\Lambda_0$, in compare to the other two
methods.

\section{ Example}

We apply our method to the breast cancer study mentioned in the
introduction. Various risk factors were measured on probands and
their relatives. For illustrative purposes we consider age at
first full-term pregnancy with the relatives of the probands being
the mothers. The following analysis is based on 437 breast cancer
case probands matched with 437 control probands and a total of 874
mothers. The number of mothers who had breast cancer was 70 among
the case families and 35 among the control families. The number of
women whose first live birth occurred before age 20 was 142 among
the probands and 181 among the mothers. In the following analysis,
the gamma frailty model is used with expectation 1 and variance
$\theta$. Three estimation procedures are considered: our proposed
method, the Hsu et al.\ (2004) method, and the modified
Shih-Chatterjee (2002) method. For our proposed method, the
two-stage procedure for $\Lambda_0$ is used since the age range of
the mothers with breast cancer was 20-76 and of the age range of
the probands was 22-44. Table 6 presents the regression
coefficient parameter estimate $\hat{\beta}$, the dependency
parameter estimate, $\hat{\theta}$, and $\hat{\Lambda}_0$ at ages
40, 50, 60 and 70 years old, along with their respective bootstrap
standard errors. The proposed approach and that of Shih and
Chatterjee yielded similar dependency estimates with the proposed
approach being moderately more efficient. Hsu et al.'s approach
gave a slightly lower dependence estimate. The regression
coefficient estimates of Hsu et al.\ and that of Shih and
Chatterjee are similar, with the latter being slightly more
efficient. The proposed approach yielded a slightly lower
covariate effect. The cumulative baseline hazard estimates are
similar under the three estimation techniques. The results, based
on the three methods, imply that women who had their first
full-term pregnancy before age 20 have a reduced risk of
developing breast cancer, supporting the observation of breast
cancer risk reduced by early first full-term pregnancy (e.g.\
Coditz et al., 1996; among others). The estimates of the
dependency parameter imply that after adjusting for the first
full-term pregnancy, there remains a significant dependency
between the ages of onset for mothers and daughters with cross
ratio ($1+\theta$) close to 2.

\section{ Discussion}
In this work we have presented a new estimator for case-control
family study survival data under a frailty model, allowing an arbitrary
frailty distribution with finite moments.
Rigorous large sample theory has been provided.
Simulation results
under the popular gamma frailty model indicate that the proposed
procedure provides estimates with minimal bias and confidence
intervals with the appropriate coverage rate. Moreover, our
estimators were seen to be essentially identical in efficiency
to estimators based on the more complex approach of
of Shih and Chatterjee (2002).

Rigorous large sample theory has been provided for unrestricted
sampling of probands. For restricted sampling, the asymptotic
theory could be worked out largely following the arguments for the
two-stage estimator but the algebra becomes very complicated.  It
is beyond the scope of the current paper and will be presented in
a future communication.

\section{ Appendix: Asymptotic theory}

This appendix presents the technical conditions we assume for the asymptotic results
and the proofs of these results. The development is patterned after Zucker (2005)
and Zucker et al.\ (2006), but considerable extension of the arguments is required.
In the presentation below, we focus on the added arguments needed for the present
setting, and refer back to Zucker (2005) and Zucker et al.\ (2006) for the other
segments of the development.

\subsection{ Assumptions and background}
In deriving the asymptotic properties of $\hat{\bgamma}$, we make
a number of assumptions. Several of these assumptions have already
been presented in the main text. Below we list the additional
assumptions.

\begin{enumerate}
\item
There is a finite maximum follow-up time $\tau>0$, with
$\E[\sum_{j=1}^{m_i}Y_{ij}(\tau)]=y^*>0$ for all $i$.
\item
The frailty random variable $\omega_i$ has finite moments up to order
$(m+2)$.
\item
There exist $b>0$ and $C>0$ such that
$$
\lim_{w \rightarrow 0} w^{-(b-1)} f(w) = C.
$$
\item
The baseline hazard function $\lambda_0^\circ(t)$ is bounded over $[0,\tau]$
by some fixed (but not necessarily known) constant $\lambda_{max}$.
\item
The function $f^\prime(w;\theta) = (d/d\theta) f(w;\theta)$ is absolutely
integrable.
\item
For any given family, there is a positive probability of at least two
failures.
\item
Defining $\pi(s) = \E [I(T_{i0} < s) \sum_{j=1}^{m_i} Y_{ij}(s)]$, we have
\begin{equation}
\xi_r(u) \equiv \int_0^u \frac{{\lambda}_0^\circ(t)}{\pi(s)^r} \, ds < \infty
\quad \mbox{for all } u \in [0,\tau] \mbox{ and } r = 1, 2, 3.
\label{icon}
\end{equation}
This assumption is needed in the analysis of the first-stage estimator.
For $r=1$, it parallels Assumption (5.4) of Keiding and Gill (1990),
\item
The matrix $[(\partial/\partial \bgamma)
\bU(\bgamma,\hat{\Lambda}_0(\cdot,\bgamma))]|_{\bgams=\bgams^\circ}$
is invertible with probability going to 1 as $n \rightarrow
\infty$.
\end{enumerate}

\subsection{Technical Preliminaries}

With $\psi^*(r,h)$ as in (\ref{psi}), we define
$\psi_{min}^*(h) = \min_{0 \leq r \leq m} \psi^*(r,h)$ and
$\psi_{max}^*(h)$ \linebreak[4]
= $\max_{0 \leq r \leq m} \psi^*(r,h)$.
In (\ref{psi}), the numerator and denominator
are bounded above since $\omega_i$ is assumed to have finite
$(m+2)$-th moment. Also, since $\omega_i$ is nondegenerate,
the numerator and denominator are strictly positive.
Thus $\psi_{max}^*(h)$
is finite and $\psi_{min}^*(h)$ is strictly positive.
We present below two lemmas. The first lemma,
which can be proved by elementary
calculus, is taken from Zucker et al.\ (2006).
The second lemma parallels Lemma 3 of Zucker et al.\ (2006).

\vspace*{0.5em}
\noindent {\bf Lemma 1:} The function $\psi^*(r,h)$ is decreasing in $h$.
Hence for all $\bgamma \in \cG$ and all $t$,
$\psi_i(\bgamma,\Lambda,t) \leq \psi_{max}^*(0)$ and
$\psi_i(\bgamma,\Lambda,t) \geq \psi_{min}^*(m \nu \Lambda(t))$.
In addition, there exist $B>0$ and $\bar{h}>0$ such that,
for all $h \geq \bar{h}$,
$\psi_{min}^*(h) \geq B h^{-1}$.

\vspace*{0.5cm}
\noindent {\bf Lemma 2:} For any $\eps>0$, we have $\sup_{s \in [\eps,\tau]}|\tilde{\Lambda}_0(s,\bgamma^\circ)
-\tilde{\Lambda}_0(s-,\bgamma^\circ)| \rightarrow 0 $ as $n \rightarrow \infty$.

\subsection{Consistency}

As indicated in Sec.~4, the consistency proof proceeds in several stages.

\noindent
{\bf Claim A1:} $\tilde{\Lambda}_0(t,\bgamma)$ converges in probability to
some function $\Lambda_0^*(t,\bgamma)$ uniformly in $t$ and $\bgamma$. The
function $\Lambda_0^*(t,\bgamma)$ satisfies $\Lambda_0^*(t,\bgamma^\circ)=\Lambda_0^\circ(t)$.

{\it Remark}: We give here an in pr.\ consistency result, rather than an a.s.\ result
as in Zucker et al.\ (2006). The reason will be explained in the course of the proof.

\noindent
{\bf Proof:} We can write $\tilde{\Lambda}_0(t,\bgamma)$ as
\begin{equation}
\tilde{\Lambda}_0(t,\bgamma)
= \int_0^t \frac
{n^{-1} \sum_{i=1}^{n} I(T_{i0} < s) \sum_{j=1}^{m_i}dN_{ij}(s)}
{n^{-1} \sum_{i=1}^{n} I(T_{i0} < s)
\bar{\psi}_i(s-,\bgamma,\tilde{\Lambda})
\sum_{j=1}^{m_i} Y_{ij}(s)\exp(\bbeta^T {\bf Z}_{ij})} \, .
\label{lt}
\end{equation}
The proof here builds here on that of the corresponding Claim A in Zucker et al. The
main point needing attention here is the fact that, because of the indicators
$I(T_{i0}<s)$, the denominator of (\ref{lt}) tends to 0 as $s \rightarrow 0$. Special
arguments are needed to deal with this ``vanishing denominator" problem.

Define, in parallel with Zucker et al.\ (2006),
$$
\Xi_n(t,\bgamma,\Lambda)=\int_0^t
\frac{n^{-1} \sum_{i=1}^{n}
I(T_{i0} < s) \sum_{j=1}^{m_i}dN_{ij}(s)}
{n^{-1} \sum_{i=1}^{n} I(T_{i0} < s)
\bar{\psi}_i(s-,\bgamma,\Lambda)
\sum_{j=1}^{m_i} Y_{ij}(s)\exp(\bbeta^T {\bf Z}_{ij})}
$$
and
$$
\Xi(t,\bgamma,\Lambda)=\int_0^t
\frac{\E [I(T_{i0}<s) \bar{\psi}_i(s-,\bgamma,\Lambda_0^\circ)
\sum_{j=1}^{m_i} Y_{ij}(s)\exp(\bbeta^{\circ T}{\bf Z}_{ij})]}
{\E [I(T_{i0}<s)
\bar{\psi}_i(s-,\bgamma,\Lambda)
\sum_{j=1}^{m_i}Y_{ij}(s) \exp(\bbeta^T {\bf Z}_{ij})]}
\lambda_0^\circ(s) ds.
$$
By definition, $\tilde{\Lambda}_0(t,\bgamma)$ satisfies the equation
$\tilde{\Lambda}_0(t,\bgamma) = \Xi_n(t,\bgamma,\hat{\Lambda}_0(\cdot,\bgamma))$.

\it Remark\rm: In Zucker et al.\ (2006), we had a result to the effect that $\Xi_n(t,\bgamma,\Lambda)
\rightarrow \Xi(t,\bgamma,\Lambda)$ a.s.\ as $n \rightarrow \infty$, uniformly over $t \in [0,\tau]$,
$\bgamma \in \cG$, and $\Lambda$ in a certain set. We could not obtain the corresponding result here;
the argument of Aalen (1976) fails in the neighborhood of zero because of the vanishing denominator
problem. This is why we give only an in pr.\ consistency result rather than an a.s.\ result.

Again in parallel with Zucker et al.\ (2006), define
$$
q_{\bgams}(s,\Lambda)
= \frac{\E [I(T_{i0}<s)
\bar{\psi}_i(s-,\bgamma,\Lambda_0^\circ)
\sum_{j=1}^{m_i} Y_{ij}(s)\exp(\bbeta^{\circ T}{\bf Z}_{ij})]}
{\E [I(T_{i0}<s)
\bar{\psi}_i(s-,\bgamma,\Lambda)
\sum_{j=1}^{m_i} Y_{ij}(s) \exp(\bbeta^T {\bf Z}_{ij})]} \lambda_0^\circ(s).
$$
This function $q_{\bgams}(s,\Lambda)$ has the same properties as noted for the corresponding
function in Zucker et al. These properties are not interfered with by the insertion of the indicator function
$I(T_{i0}<s)$. In particular, from Lemma 1 we have
$$
\frac{\bar{\psi}_i(s-,\bgamma,\Lambda_0^\circ)}
{\bar{\psi}_i(s-,\bgamma,\Lambda)}
\leq \frac{\psi_{max}^*(0)}{\psi_{min}^*(h_{max})},
$$
Pulling this bound outside of the expectation, we get the a bound on
$q_{\bgams}(s,\Lambda)$ analogous to that in Zucker et al.
Similarly, as in Zucker et al., the function $q_{\bgams}(s,\Lambda)$
has the following Lipschitz-like property:
$|q_{\bgams}(s,\Lambda_1)-q_{\bgams}(s,\Lambda_2)| \leq
K \sup_{0 \leq u \leq s} |\Lambda_1(u)-\Lambda_2(u)|$.
Accordingly, we find that the
equation $\Lambda(t) = \Xi(t,\bgamma,\Lambda)$ has a unique solution, which we denote
by $\Lambda_0^*(t,\bgamma)$. The claim then is that $\tilde{\Lambda}_0(t,\bgamma)$
converges in pr.\ (uniformly in $t$ and $\bgamma$) to $\Lambda_0^*(t,\bgamma)$.

We now define, for any $\epsilon>0$, the quantities
$$
\Xi_n(t,\bgamma,\Lambda,\epsilon)=\int_\epsilon^t
\frac{n^{-1} \sum_{i=1}^{n}
I(T_{i0} < s) \sum_{j=1}^{m_i}dN_{ij}(s)}
{n^{-1} \sum_{i=1}^{n} I(T_{i0} < s)
\bar{\psi}_i(s-,\bgamma,\Lambda)
\sum_{j=1}^{m_i} Y_{ij}(s)\exp(\bbeta^T {\bf Z}_{ij})}
$$
and
$$
\Xi(t,\bgamma,\Lambda,\epsilon)=\int_\epsilon^t
\frac{\E [I(T_{i0}<s) \bar{\psi}_i(s-,\bgamma,\Lambda_0^\circ)
\sum_{j=1}^{m_i} Y_{ij}(s)\exp(\bbeta^{\circ T}{\bf Z}_{ij})]}
{\E [I(T_{i0}<s)
\bar{\psi}_i(s-,\bgamma,\Lambda)
\sum_{j=1}^{m_i}Y_{ij}(s) \exp(\bbeta^T {\bf Z}_{ij})]}
\lambda_0^\circ(s) ds.
$$
We next define $\tilde{\Lambda}_0(t,\bgamma,\epsilon)$ to be the solution of the equation
$\tilde{\Lambda}_0(t,\bgamma,\epsilon) = \Xi_n(t,\bgamma,\tilde{\Lambda}_0(\cdot,\bgamma),\epsilon)$,
starting from $\tilde{\Lambda}_0(\epsilon,\bgamma,\epsilon)=0$. We extend the definition of
$\tilde{\Lambda}_0(t,\bgamma,\epsilon)$ by setting it equal to 0 for $t<\epsilon$. Similarly,
we define $\Lambda_0^*(t,\bgamma,\epsilon)$ to be the solution of the equation
$\Lambda_0(t,\bgamma,\epsilon) = \Xi(t,\bgamma,\Lambda_0(\cdot,\bgamma),\epsilon)$,
starting from $\Lambda_0(\epsilon,\bgamma,\epsilon)=0$, and extend the definition
by setting $\Lambda_0^*(t,\bgamma,\epsilon)$ equal to 0 for $t<\epsilon$.

For $t \in [\epsilon,\tau]$, the difference between $\Lambda_0^*(t,\bgamma)$ and
$\Lambda_0^*(t,\bgamma,\epsilon)$ is as follows: $\Lambda_0^*(t,\bgamma)$ is the solution
to $\Lambda_0(t,\bgamma,\epsilon) = \Xi(t,\bgamma,\Lambda_0(\cdot,\bgamma),\epsilon)$,
starting from $\Lambda_0(\epsilon,\bgamma,\epsilon)=\Lambda_0(\epsilon,\bgamma)$, whereas
$\Lambda_0^*(t,\bgamma,\epsilon)$ is the solution to $\Lambda_0(t,\bgamma,\epsilon) =
\Xi(t,\bgamma,\Lambda_0(\cdot,\bgamma),\epsilon)$, starting from $\Lambda_0(\epsilon,
\bgamma,\epsilon)$ $=0$. Hence, by an induction argument similar to that in the proof of
Hartman (1973, Thm.\ 1.1), we find that
$$
|\Lambda_0^*(t,\bgamma,\epsilon) - \Lambda_0^*(t,\bgamma)|
\leq e^K \Lambda_0(\epsilon,\bgamma),
$$
where $K$ is the Lipschitz constant for $q_{\bgams}(s,\Lambda)$. We thus have
\begin{equation}
\sup_{\bgams \in \cG, \, t \in [0,\tau]} |\Lambda_0^*(t,\bgamma,\epsilon) - \Lambda_0^*(t,\bgamma)|
\rightarrow 0 \quad \mbox{as } \epsilon \rightarrow 0.
\label{oops}
\end{equation}

Now, for any given $\epsilon>0$, there is no vanishing denominator problem on the
interval $[\epsilon,\tau]$. Hence, the argument in Zucker et al.\ (2006) goes through
as is, and we get the following result: for any $\epsilon > 0$,
\begin{equation}
\sup_{\bgams \in \cG, \, t \in [\epsilon,\tau]}
|\tilde{\Lambda}_0(t,\bgamma,\epsilon) -\Lambda_0^*(t,\bgamma,\epsilon)|
\rightarrow 0 \quad \mbox{a.s.\ as } n \rightarrow \infty.
\end{equation}
In fact, in the supremum above, we can replace $[\epsilon,\tau]$ by $[0,\tau]$,
since by definition $\tilde{\Lambda}_0(t,\bgamma,\epsilon)=\Lambda_0(t,\bgamma,\epsilon)=0$ for
$t < \epsilon$.

The above a.s.\ result immediately yields the corresponding in pr.\ result:
\begin{equation}
\sup_{\bgams \in \cG, \, t \in [0,\tau]}
|\tilde{\Lambda}_0(t,\bgamma,\epsilon) -\Lambda_0^*(t,\bgamma,\epsilon)|
\rightarrow 0 \quad \mbox{in pr.\ as } n \rightarrow \infty.
\label{ce}
\end{equation}

Our aim now is to show that
\begin{equation}
\sup_{\bgams \in \cG, \, t \in [0,\tau]}
|\tilde{\Lambda}_0(t,\bgamma) -\Lambda_0^*(t,\bgamma)|
\rightarrow 0 \quad \mbox{in pr.\ as } n \rightarrow \infty.
\end{equation}
That is, we want to show the following: for any $\rho,\delta>0$, there exists $n^*(\rho,\delta)$
large enough such that
$$
\Pr (\sup_{\bgams \in \cG, \, t \in [0,\tau]}
|\tilde{\Lambda}_0(t,\bgamma) -\Lambda_0^*(t,\bgamma)| > \rho) \leq \delta
$$
for all $n \geq n^*(\rho,\delta)$.

Let $\rho$ and $\delta$ be given. By (\ref{oops}), we can find $\epsilon>0$ small
enough such that
$$
\sup_{\bgams \in \cG, \, t \in [0,\tau]} |\Lambda_0(t,\bgamma,\epsilon) - \Lambda_0^*(t,\bgamma)|
\leq \frac{\rho}{3} \, .
$$
Further, for this fixed $\epsilon$, the result (\ref{ce}) implies that there exists
$\tilde{n}$ such that
$$
\Pr \left( \sup_{\bgams \in \cG, \, t \in [0,\tau]}
|\tilde{\Lambda}_0(t,\bgamma,\epsilon) -\Lambda_0^*(t,\bgamma,\epsilon)|
> \frac{\rho}{3} \right) \leq \frac{\delta}{2}
$$
for all $n \geq \tilde{n}$.

Now

$|\tilde{\Lambda}_0(t,\bgamma) -\Lambda_0^*(t,\bgamma)|$
\begin{equation}
\leq
|\tilde{\Lambda}_0(t,\bgamma) -\tilde{\Lambda}_0(t,\bgamma,\epsilon)|
+ |\tilde{\Lambda}_0(t,\bgamma,\epsilon) -\Lambda_0^*(t,\bgamma,\epsilon)|
+ |\Lambda_0^*(t,\bgamma,\epsilon) -\Lambda_0^*(t,\bgamma)|
\label{string}
\end{equation}
The developments just above imply that, for $n \geq \tilde{n}$, the supremum over
$\bgamma \in \cG$ and $t \in [0,\tau]$ of the sum of the last two terms is bounded
by $\frac{2}{3}\rho$ with probability at least $1-\half \delta$. It remains to deal
with the first term.

Define
\begin{eqnarray*}
C_1(s) &=&
\frac{1}{n} \sum_{i=1}^{n} I(T_{i0} < s)
\bar{\psi}_i(s-,\bgamma,\tilde{\Lambda}(\cdot,\bgamma))
\sum_{j=1}^{m_i} Y_{ij}(s)\exp(\bbeta^T {\bf Z}_{ij}), \\
C_2(s) &=&
\frac{1}{n} \sum_{i=1}^{n} I(T_{i0} < s)
\bar{\psi}_i(s-,\bgamma,\tilde{\Lambda}(\cdot,\bgamma,\epsilon))
\sum_{j=1}^{m_i} Y_{ij}(s)\exp(\bbeta^T {\bf Z}_{ij}).
\end{eqnarray*}
We can then write

%\newpage

\begin{equation}
\tilde{\Lambda}_0(t,\bgamma) -\tilde{\Lambda}_0(t,\bgamma,\epsilon)
= \tilde{\Lambda}_0(t \wedge \epsilon,\bgamma) + A(t,\epsilon),
\label{wow}
\end{equation}
where
$$
A(t,\epsilon) =  \int_{t \wedge \epsilon}^t
[C_1(s)^{-1} - C_2(s)^{-1}] \left[ \frac{1}{n} \sum_{i=1}^{n} I(T_{i0} < s)
\sum_{j=1}^{m_i} dN_{ij}(s) \right].
$$
We deal with the two terms on the right side of (\ref{wow}) in turn. In what follows,
we let $R$ denote a ``generic" constant which may vary from one appearance to another,
but does not depend on the unknown parameters or $\epsilon$.

Denote
$$
\Pi(s) =
\frac{1}{n} \sum_{i=1}^{n} I(T_{i0} < s) \sum_{j=1}^{m_i} Y_{ij}(s)
$$
and recall the definition $\pi(s) = \E [I(T_{i0} < s) \sum_{j=1}^{m_i} Y_{ij}(s)]$.
Also recall
$$
\tilde{\Lambda}(t,\bgamma)
= \int_0^t
\frac
{\frac{1}{n}
\sum_{i=1}^{n} I(T_{i0} < s)
\sum_{j=1}^{m_i}dN_{ij}(s)}
{\frac{1}{n} \sum_{i=1}^{n} I(T_{i0} < s)
\bar{\psi}_i(s-,\bgamma,\tilde{\Lambda}(\cdot,\bgamma))
\sum_{j=1}^{m_i} Y_{ij}(s)\exp(\bbeta^T {\bf Z}_{ij})}
$$
It is clear that $\tilde{\Lambda}(t,\bgamma) \leq R \Upsilon(t,\bgamma)$, where
$$
\Upsilon(t,\bgamma)
= \int_0^t
\frac
{\frac{1}{n}
\sum_{i=1}^{n} I(T_{i0} < s) \sum_{j=1}^{m_i}dN_{ij}(s)}
{\frac{1}{n} \sum_{i=1}^{n} I(T_{i0} < s) \sum_{j=1}^{m_i} Y_{ij}(s)}
= \int_0^t
\frac
{\frac{1}{n}
\sum_{i=1}^{n} I(T_{i0} < s) \sum_{j=1}^{m_i}dN_{ij}(s)}
{\Pi(s)}.
$$
We can write
$$
\Upsilon(t,\bgamma) = \int_0^t
\Pi(s)^{-1}
\left[ \frac{1}{n}
\sum_{i=1}^{n} I(T_{i0} < s)
\sum_{j=1}^{m_i} Y_{ij}(s)\exp({\bbeta^\circ}^T {\bf Z}_{ij}) \right]
\lambda_0^\circ(s) ds
$$
$$
+
\int_0^t
\Pi(s)^{-1}
\left[ \frac{1}{n}
\sum_{i=1}^{n} I(T_{i0} < s)
\sum_{j=1}^{m_i} dM_{ij}(s) \right],
$$
where $M_{ij}$ is the martingale process corresponding to $N_{ij}$:
\begin{equation}
M_{ij}(t)=N_{ij}(t)-\int_0^t
\lambda_0(u) \exp(\bbeta^{\circ T} {\bf Z}_{ij}) Y_{ij}(u) \psi_{i}(\bgamma^\circ,\Lambda_0^\circ,u-) du.
\label{mart}
\end{equation}

The first term is
clearly bounded by $R \Lambda_0^\circ(t)$. Thus, denoting the second term by $M^*(t)$,
we have
\begin{equation}
\tilde{\Lambda}(t,\bgamma)
\leq R [\Lambda_0^\circ(t) + \sup_{u \in [0,\tau]} |M^*(u)|]
\label{bou}
\end{equation}

We next examine $A(t,\epsilon)$. We can restrict to $t \geq \epsilon$,
since $A(t,\epsilon)=0$ for $t < \epsilon$.
Denote $\Delta(t) = \tilde{\Lambda}_0(t,\bgamma) -
\tilde{\Lambda}_0(t,\bgamma,\epsilon)$.
Bearing in mind the Lipschitz property
of $\bar{\psi}$, we find that
$$
|A(t,\epsilon)|
\leq R \int_{\epsilon}^t
|\Delta(s-)| d\Upsilon(s).
$$
Note that, for $t \geq \epsilon$, $dA(t,\epsilon)=d\Delta(t)$. Thus,
a simple induction and some additional simple manipulations lead to the following,
where we employ the symbol $\cP$ to denote product integral and use the fact that
$\Delta(\epsilon)=\tilde{\Lambda}(\epsilon,\bgamma)$:
$$
|A(t,\epsilon)|
\leq
|\Delta(\epsilon)| \cP_\epsilon^t (1 + R d\Upsilon(s))
\leq
|\Delta(\epsilon)| \exp(R [\Upsilon(t)-\Upsilon(\epsilon)])
\leq
|\tilde{\Lambda}(\epsilon,\bgamma)| \exp(R \Upsilon(\tau))
$$
In view of the analysis above of $\Upsilon(t)$, we get
\begin{equation}
|A(t,\epsilon)|
\leq
|\tilde{\Lambda}(\epsilon,\bgamma)|
\exp (R [\Lambda_0^\circ(\tau) + \sup_{u \in [0,\tau]} |M^*(u)|]).
\label{boutwo}
\end{equation}

Putting (\ref{wow}), (\ref{bou}), and (\ref{boutwo}) together, we get

$|\tilde{\Lambda}_0(t,\bgamma) -\tilde{\Lambda}_0(t,\bgamma,\epsilon)|$
\begin{equation}
\leq
R_1 [\Lambda_0^\circ(\epsilon) + \sup_{u \in [0,\tau]} |M^*(u)|]
(1 + \exp (R_2[1 + \sup_{u \in [0,\tau]} |M^*(u)|]).
\label{almost}
\end{equation}
for suitable absolute constants $R_1$ and $R_2$.

The last main step is to analyze the martingale process
$$
M^*(u)
= \int_0^u
\Pi(s)^{-1}
\left[ \frac{1}{n}
\sum_{i=1}^{n} I(T_{i0} < s)
\sum_{j=1}^{m_i} dM_{ij}(s) \right].
$$
Our argument is pattered after the argument given by
Keiding and Gill (1990, p.\ 595).

By Lenglart's and Markov's inequalities, we have, for any positive $\kappa$ and $\eta$
and any $c \in [0,\tau]$,
$$
\Pr(\sqrt{n} \, \sup_{u \in [0, c]} |M^*(t)| > \kappa)
\leq \eta + \Pr(n \langle M^* \rangle (c) > \eta/\kappa^2)
\leq \eta +  \frac{\kappa^2}{\eta} \E [n \langle M^* \rangle (c)].
$$
Define $J(s) = I(\Pi(s)>0)$. Then
$$
n \langle M^* \rangle (c)
= \int_0^c
\Pi(s)^{-2}
\left[ \frac{1}{n}
\sum_{i=1}^{n} I(T_{i0} < s)
\sum_{j=1}^{m_i} Y_{ij}(s)\exp({\bbeta^\circ}^T {\bf Z}_{ij}) \right] \lambda_0^\circ(s) ds
$$
$$
\leq
R \int_0^c \frac{nJ(s)}{n\Pi(s)} \lambda_0^\circ(s) ds
\leq
R \int_0^c \frac{n+1}{n\Pi(s)+1} \lambda_0^\circ(s) ds.
$$
As in Keiding and Gill,
$$
\E \left[ \frac{n+1}{n\Pi(s)+1} \right] \leq \frac{1}{\pi(s)} \, .
$$
Hence $\E [n \langle M^* \rangle (c) ] \leq R \xi_1(c)$, where $\xi_1$ is as in Assumption 7.
We thus get
\begin{equation}
\Pr(\sup_{u \in [0, c]} |M^*(u)| > \kappa n^{-\half})
\leq \eta + R \kappa^2 \eta^{-1} \xi_1(c).
\label{mbou}
\end{equation}
Now, the main quantities in the bound in (\ref{almost}) are $\Lambda_0^\circ(\epsilon)$
and $\sup_{u \in [0, \tau]} |M^*(u)|$. By decreasing $\epsilon$ if necessary, we can make
$\Lambda_0^\circ(\epsilon)$ as small as we need. The behavior of $\sup_{u \in [0, \tau]}
|M^*(u)|$ is characterized by (\ref{mbou}). We see that by decreasing $\epsilon$ if necessary
and choosing $\eta$ appropriately, we can guarantee that the probability that the right side of
(\ref{almost}) is less than $\frac{1}{3}\rho$ will be at least $1-\half \delta$ for all $n$
sufficiently large. With this, we have taken care of the first term of (\ref{string}). The
desired convergence has thus been established. It is easy to see that $\Lambda_0^*(t,\bgamma^\circ)
=\Lambda_0^\circ(t)$.

\noindent
{\bf Claim A2:} $\hat{\Lambda}_0(t,\bgamma)$ converges in probability to
some function $\Lambda_0(t,\bgamma)$ uniformly in $t$ and $\bgamma$.
The function $\Lambda_0(t,\bgamma)$ satisfies $\Lambda_0(t,\bgamma^\circ)=
\Lambda_0^\circ(t)$.

\noindent
{\bf Proof:}
We can write $\hat{\Lambda}_0(t,\bgamma)$ as
$$
\hspace*{-3em} \hat{\Lambda}_0(t,\bgamma)
= \int_0^t \frac
{n^{-1} \sum_{i=1}^{n} \sum_{j=1}^{m_i}dN_{ij}(s)}
{n^{-1} \sum_{i=1}^{n} \tilde{\psi}_i(s-,\bgamma)
\sum_{j=1}^{m_i} Y_{ij}(s)\exp(\bbeta^T {\bf Z}_{ij})}.
$$
In view of Claim A1 above, up to a uniform error of $o_P(1)$ we can
replace all instances of $\tilde{\Lambda}_0(u,\bgamma)$ in the definition
of $\tilde{\psi}_i(s-,\bgamma)$ by $\Lambda_0^*(u,\bgamma)$. The desired
result then can be obtained using the argument used to prove Claim A
of Zucker et al.\ (2006).

\noindent
{\bf Claim A3:} We have
\begin{eqnarray*}
\sup_{s \in [0,\tau], \bgams \in \cG}|\tilde{\Lambda}_0(s,\bgamma^\circ)
-\tilde{\Lambda}_0(s-,\bgamma^\circ)| & \CP & 0 \mbox{ as } n \rightarrow \infty, \\
\sup_{s \in [0,\tau], \bgams \in \cG}|\hat{\Lambda}_0(s,\bgamma^\circ)
-\hat{\Lambda}_0(s-,\bgamma^\circ)| & \CP & 0 \mbox{ as } n \rightarrow \infty.
\end{eqnarray*}

\noindent
{\bf Proof:}
These results follow from Claims A1 and A2 and the fact that $\Lambda_0^*(t,\bgamma)$
and $\Lambda_0(t,\bgamma)$ are continuous.

\noindent
{\bf Claim B:}
${\bf U}(\bgamma,\hat{\Lambda}_0(\cdot,\bgamma))$ converges in probability
uniformly in $t$ and $\bgamma$ to a limit ${\bf u} (\bgamma,\Lambda_0(\cdot,\bgamma))$.

\noindent {\bf Proof:} As in Claim B of Zucker et al.\ (2006).

\noindent
{\bf Claim C:} There exists a
unique consistent (in pr.) root to ${\bf U}(\hat{\bgamma},\hat{\Lambda}_0(\cdot,\hat{\bgamma}))=\bzro$.

\noindent {\bf Proof:} By appeal to Foutz's (1977) theorem, as in Claim C of Zucker et al.\ (2006).

\subsection{A workable representation of $\hat{\Lambda}_0(t)-\Lambda_0^\circ(t)$}

In order to develop our asymptotic normality result, we need a workable representation
of $\hat{\Lambda}_0(t)-\Lambda_0^\circ(t)$. The first step is to develop a suitable
representation of $\tilde{\Lambda}_0(t)-\Lambda_0^\circ(t)$. Then, building on this,
we develop our representation of $\hat{\Lambda}_0(t)-\Lambda_0^\circ(t)$.

\vspace*{1em}
\noindent
{\bf 9.4.1 Representation of $\tilde{\Lambda}_0(t)-\Lambda_0^\circ(t)$}

Our starting point is the following simple lemma.

%\newpage

\noindent
{\bf Lemma:}
Let $\mathcal{R}_n(t)$ and $\mathcal{S}_n(t)$ be stochastic processes, and let
$A_n(t,\eps)$ and $B_n(t,\eps)$ be quantities that are bounded in probability
uniformly in $t$ and $\eps$. Define
\begin{eqnarray*}
\mathcal{R}_n(t,\eps)&=&\mathcal{R}_n(t)-A_n(t,\eps)\mathcal{R}_n(\eps), \\
\mathcal{S}_n(t,\eps)&=&B_n(t,\eps)[\mathcal{S}_n(t)-\mathcal{S}_n(\eps)].
\end{eqnarray*}
Suppose that:
\begin{enumerate}
\item
$\sup_{t \in [\eps,\tau]} \sqrt{n} |\mathcal{R}_n(t,\eps) - \mathcal{S}_n(t,\eps)|
\CP 0$ as $n \rightarrow \infty$ for any fixed $\eps>0$.
\item
$\lim_{\eps \downarrow 0} \limsup_{n \rightarrow \infty}
\Pr(\sup_{t \in [0,\eps]} \sqrt{n} |\mathcal{R}_n(t)| > \delta) = 0$ for all $\delta > 0$.
\item
$\lim_{\eps \downarrow 0} \limsup_{n \rightarrow \infty}
\Pr(\sup_{t \in [0,\eps]} \sqrt{n}  |\mathcal{S}_n(t)| > \delta) = 0$ for all $\delta > 0$.
\item
$B_n(t,\eps) \rightarrow B_n(t,0)$ uniformly in $t$ as $\eps \rightarrow 0$ with probability
converging to one as $n \rightarrow \infty$.
\end{enumerate}
Then $\sup_{t \in [0,\tau]} \sqrt{n} |\mathcal{R}_n(t) - B_n(t,0)\mathcal{S}_n(t,0)| \CP 0$.

\vspace*{1em}
We apply this lemma with $\mathcal{R}_n(t) = \sqrt{n}[\tilde{\Lambda}_0(t)-\Lambda_0^\circ(t)]$.
We have to check the four conditions enumerated in the lemma.

\bsh \underline{Condition 1} \esh

Arguments along the lines of Zucker et al.\ (2006) yield the result of Condition 1, with
\begin{eqnarray}
\mathcal{S}_n(t) &=&
\int_0^t
\frac{\tilde{p}(s-,\eps)}
{\tilde{\mathcal{Y}}(s,\Lambda_0^\circ)}
\left[ \frac{1}{n} \sum_{i=1}^{2n}\sum_{j=1}^{m_i}I(T_{i0}<s) dM_{ij}(s) \right]
, \label{eq:martin} \\
A_n(t,\eps) & = & B_n(t,\eps) = \tilde{p}(t,\eps)^{-1}, \nonumber
\end{eqnarray}
where
\begin{equation}
\tilde{p}(t,\eps)=\prod_{s \in [\eps,t]} \left[
 1+n^{-1}\sum_{i=1}^{2n} \sum_{j=0}^{m_i}
 \Omega_{ij}(s,t)  d\tilde{N}_{ij}(s) + n^{-1} \Omega^*(s)I(T_{i0}<s)\delta_{ij}
 \right].
\end{equation}
Here
$$
 \tilde{\mathcal{Y}}(s,\Lambda) = \frac{1}{n}  \sum_{i=1}^{2n}
 I(T_{i0}<s)\psi_i(\bgamma^\circ,\Lambda,s)R_{i.}(s)
$$
with $ R_{i.}(s)=\sum_{j=1}^{m_i} Y_{ij}(s)
 \exp(\bbeta^{\circ T} {\bf Z}_{ij})$,
$$
\Omega^*(s) = \frac{1}{n}\sum_{k=1}^{2n}\frac{R_{k.}(s)
\eta_{1k}(s)I(T_{k0}<s)}{\{\tilde{\mathcal{Y}}(s,\Lambda_0^\circ)\}^{2}}
\sum_{l=1}^{m_k}I(T_{kl}>s)\exp({\bbeta^{\circ T} {\bf Z}_{kl}}),
$$
$$
\Omega_{i0}(s,t)=n^{-1}\int_s^t \frac{R_{i.}(u) \eta_{1i}(u)
\exp(\bbeta^{\circ T} {\bf Z}_{i0})}{\{\tilde{\mathcal{Y}}(u,\Lambda_0^\circ)\}^{2}}
\sum_{k=1}^{2n}\sum_{l=1}^{m_k}d N_{kl}(u),
$$
and for $j \geq 1$
$$
{\Omega_{ij}(s,t)=n^{-1}\int_s^t\frac{I(T_{i0}<u)
R_{i.}(u) \eta_{1i}(u)\exp({\bbeta^{\circ T} {\bf
Z}_{ij}})}{\{\tilde{\mathcal{Y}}(u,\Lambda_0^\circ)\}^{2}}
\sum_{k=1}^{2n}\sum_{l=1}^{m_k}I(T_{k0}<u)d N_{kl}(u)}.
$$
In the above, $\eta_{1i}(s)$ is defined as
\begin{eqnarray*}
\eta_{1i}(s)=\frac{\phi_{3i}(\bgamma^\circ,{\Lambda}_0^\circ,s)}
{\phi_{1i}(\bgamma^\circ,{\Lambda}_0^\circ,s)}-\left\{
\frac{\phi_{2i}(\bgamma^\circ,{\Lambda}_0^\circ,s)}{\phi_{1i}(\bgamma^\circ,{\Lambda}_0^\circ,s)}
\right\}^2.
\end{eqnarray*}
In Sec.~9.3.2 below, we present in detail a similar argument for $\hat{\Lambda}_0(t)-\Lambda_0^\circ(t)$.

Appealing to Assumption 7 and using arguments similar to those used in the consistency proof, we find
that the $\Omega$ quantities defined above converge in probability uniformly in $s$ and $t$, so that
$\tilde{p}(t,\eps)$ converges in probability to a deterministic limit uniformly in $t$ and $\eps$.

\bsh \underline{Condition 2, 3, and 4} \esh

In regard to Condition 2, we have
$$
\tilde{\Lambda}_0(t,\bgamma) - \Lambda_0^\circ(t)
= \Delta_1(t) + \Delta_2(t),
$$
where
\begin{eqnarray*}
\Delta_1(t) & = & \int_0^t [\Gamma(s,\bgamma)-1] \lambda_0^\circ(s) ds, \\
\Delta_2(t) & = & \int_0^t \frac
{n^{-1} \sum_{i=1}^{n} I(T_{i0} < s) \sum_{j=1}^{m_i}dM_{ij}(s)}
{n^{-1} \sum_{i=1}^{n} I(T_{i0} < s)
\bar{\psi}_i(s-,\bgamma,\tilde{\Lambda})
\sum_{j=1}^{m_i} Y_{ij}(s)\exp(\bbeta^T {\bf Z}_{ij})},
%\label{lt}
\end{eqnarray*}
where
$$
\Gamma(s,\bgamma)
= \frac
{n^{-1} \sum_{i=1}^{n} I(T_{i0} < s)
\psi_i(s-,\bgamma^\circ,\Lambda_0^\circ)
\sum_{j=1}^{m_i} Y_{ij}(s)\exp({\bbeta^\circ}^T {\bf Z}_{ij})}
{n^{-1} \sum_{i=1}^{n} I(T_{i0} < s)
\bar{\psi}_i(s-,\bgamma,\tilde{\Lambda})
\sum_{j=1}^{m_i} Y_{ij}(s)\exp(\bbeta^T {\bf Z}_{ij})}
$$
and $M_{ij}(t)$ is defined as in (\ref{mart}). We will deal with $\Delta_1(t)$
and $\Delta_2(t)$ in turn, starting with $\Delta_2(t)$. In the development below,
$R$ denotes a ``generic" absolute constant.

The quadratic variation process of $\Delta_2(t)$ is given by
$$
\langle \Delta_2 \rangle (t)
= \int_0^t \left[ \frac
{n^{-1} \sum_{i=1}^{n} I(T_{i0} < s)
\psi_i(s-,\bgamma^\circ,\Lambda_0^\circ)
\sum_{j=1}^{m_i} Y_{ij}(s)\exp({\bbeta^\circ}^T {\bf Z}_{ij})}
{[n^{-1} \sum_{i=1}^{n} I(T_{i0} < s)
\bar{\psi}_i(s-,\bgamma,\tilde{\Lambda})
\sum_{j=1}^{m_i} Y_{ij}(s)\exp(\bbeta^T {\bf Z}_{ij})]^2} \right]
\lambda_0^\circ(s) ds.
$$
By arguments similar to those used in connection with $M^*(t)$ in the proof of
Claim A1, we find that $\E[n \langle \Delta_2 \rangle (t)] \leq R \xi_1(t)$.
An application of Lenglart's inequality then gives
$$
\Pr(\sqrt{n} \sup_{t \in [0, \eps]} |\Delta_2(t)| > \kappa)
\leq \eta + R \kappa^2 \eta^{-1} \xi_1(\eps) \quad \forall \eta > 0.
$$
Assumption 7 implies that $\xi_1(\eps) \downarrow 0$ as $\eps \downarrow 0$,
and this takes care of $\Delta_2(t)$.

We now turn to $\Delta_1(t)$. As before, denote $J(s) = I(\Pi(s)>0)$.
We can write
$$
\Delta_1(t) = \Delta_{1a}(t) + \Delta_{1b}(t)
$$
with
$$
\Delta_{1a}(t) = \int_0^t [\Gamma(s,\bgamma)-1]J(s) \lambda_0^\circ(s) ds
$$
and
$$
\Delta_{1b}(t) = \int_0^t [J(s)-1] \lambda_0^\circ(s) ds.
$$
The term $\Delta_{1b}(t)$ can be shown to be uniformly $O_p(n^{-\half})$
by the argument in the middle of page 595 in Keiding and Gill (1990)
As for $\Delta_{1a}(t)$, we have
$$
\Delta_{1a}(t)
\leq Rt |\tilde{\Lambda}_0(t,\bgamma) - \Lambda_0^\circ(t)|
\leq Rt |\Delta_1(t)| + Rt |\Delta_2(t)|
\leq Rt |\Delta_{1a}(t)| + Rt |\Delta_{1b}(t)| + Rt |\Delta_2(t)|.
$$
Thus, for $t$ small,
$$
|\Delta_{1a}(t)|
\leq \frac{Rt}{1-Rt} [\Delta_{1b}(t) + \Delta_2(t)],
$$
and the terms on the right hand side have already been taken care of.

The proof of Condition 3 is similar to that given above for $\Delta_2(t)$.
Condition 4 follows easily from the uniform convergence of the $\Omega$
quantities.

\vspace*{1em}
\noindent
{\bf 9.4.2 Representation of $\hat{\Lambda}_0(t)-\Lambda_0^\circ(t)$}

Let
$$
\mathcal{Y}(s,\{\tilde{\Lambda}_0,\hat{\Lambda}_0\})=\frac{1}{n}
\sum_{i=1}^{2n}\tilde{\psi}_i(\bgamma^\circ,\hat\Lambda_0,s)R_{i.}(s)
$$
and
$$
\mathcal{Y}(s,\Lambda)=\frac{1}{n}
\sum_{i=1}^{2n}{\psi}_i(\bgamma^\circ,\Lambda,s)R_{i.}(s)
$$
so that in $\tilde\psi_i(\bgamma^\circ,\hat\Lambda_0,s)$ we take
$\tilde{\Lambda}_0(T_{i0})$ if $T_{i0} \geq s$ and
$\hat{\Lambda}_0(T_{ij})$ if $T_{ij}<s$, $j \geq 0$. By Claim A3, we
have that also $\sup_{s \in [0,\tau]}
|\hat{\Lambda}_0(s,\bgamma^\circ)-\hat{\Lambda}_0(s-,\bgamma^\circ)|$
converges to zero. Thus, we obtain the following
approximation, uniformly over $t \in [0,\tau]$:
 \begin{eqnarray*}
\lefteqn{\hat{\Lambda}_0(t,\bgamma^\circ)-\Lambda_0^\circ(t)
\approx  \frac{1}{n} \int_0^t
\{{\mathcal{Y}}(s,\Lambda_0^\circ)\}^{-1} \sum_{i=1}^{2n}
 \sum_{j=1}^{m_i}  dM_{ij}(s)} \nonumber \\
 &&+ \frac{1}{n} \int_0^t \left[\{{\mathcal{Y}}(s,\{\tilde{\Lambda}_0,\hat\Lambda_0\})\}^{-1}
 - \{{\mathcal{Y}}(s,\Lambda_0^\circ)\}^{-1}\right]\sum_{i=1}^{2n}
 \sum_{j=1}^{m_i} dN_{ij}(s).
 \end{eqnarray*}
 Now let
 $$
 {\mathcal X}(s,r)=\{{{\mathcal
 Y}}(s,\Lambda_0^\circ+r\Delta^\star)\}^{-1}
 $$
with $\Delta^\star=\hat{\Lambda}_0-\Lambda_0^\circ$ or
$\tilde{\Lambda}_0-\Lambda_0^\circ$, according to the estimator
being used. Define $\dot{\mathcal X}$ and $\ddot{\mathcal X}$ as
the first and second derivative of ${\mathcal X}$ with respect to
$r$, respectively. Then, by a first order Taylor expansion of
${\mathcal X}(s,r)$ we get
\begin{eqnarray*}
 \lefteqn{\hat{\Lambda}_0(t,\bgamma^\circ)-\Lambda_0^\circ(t)  \approx
 n^{-1}\int_0^t \{{\mathcal{Y}}(s,\Lambda_0^\circ)\}^{-1} \sum_{i=1}^{2n}
 \sum_{j=1}^{m_i} dM_{ij}(s)} \\
 && \hspace*{-1.5cm}-n^{-2}\int_0^t
 \sum_{k=1}^{2n}
 \frac{R_{k.}(s) \eta_{1k}(s)}
 {\{\mathcal{Y}(s,\Lambda_0^\circ)\}^{2}}
 \sum_{l=1}^{m_k}I(T_{kl}>s)\exp(\bbeta^{\circ T} {\bf Z}_{kl})
 \{\hat{\Lambda}_0(s)-\Lambda_0^\circ(s)\}
  \sum_{i=1}^{2n}
 \sum_{j=1}^{m_i}dN_{ij}(s)  \\
 && \hspace*{-1.5cm}- n^{-2}\int_0^t
\sum_{k=1}^{2n}\frac{ R_{k.}(s)
 \eta_{1k}(s)}
 {\{\mathcal{Y}(s,\Lambda_0^\circ)\}^{2}}
 \sum_{l=1}^{m_k}I(T_{kl} \leq s)\exp(\bbeta^{\circ T} {\bf Z}_{kl})
 \{\hat{\Lambda}_0(T_{kl})-\Lambda_0^\circ(T_{kl})\}
\sum_{i=1}^{2n}
 \sum_{j=1}^{m_i} dN_{ij}(s)\nonumber\\
&& \hspace*{-1.5cm} -n^{-2} \int_0^t \sum_{k=1}^{2n}\frac{
R_{k.}(s) \eta_{1k}(s)}
 {\{\mathcal{Y}(s,\Lambda_0^\circ)\}^{2}}I(T_{k0} \geq s)  \exp(\bbeta^{\circ T} {\bf Z}_{k0})
 \{\tilde{\Lambda}_0(T_{k0})-\Lambda_0^\circ(T_{k0})\}\sum_{i=1}^{2n}
 \sum_{j=1}^{m_i}  dN_{ij}(s)\nonumber\\
 && \hspace*{-1.5cm} -n^{-2} \int_0^t \sum_{k=1}^{2n}\frac{
R_{k.}(s) \eta_{1k}(s)}
 {\{\mathcal{Y}(s,\Lambda_0^\circ)\}^{2}}I(T_{k0} < s)  \exp(\bbeta^{\circ T} {\bf Z}_{k0})
 \{\hat{\Lambda}_0(T_{k0})-\Lambda_0^\circ(T_{k0})\}\sum_{i=1}^{2n}
 \sum_{j=1}^{m_i}  dN_{ij}(s).\nonumber
\end{eqnarray*}
The justification for ignoring the remainder term in the Taylor
expansion is as in the parallel argument in Zucker et al.\ (2006).

The second, third and fifth terms of the above equation can be
written, by interchanging the order of integration, as
$$
-n^{-1}\int_0^t\{\hat{\Lambda}_0(s)-\Lambda_0^\circ(s)\}\sum_{i=1}^{2n}\sum_{j=0}^{m_i}
\Upsilon_{ij}(s,t)d\tilde N_{ij}(s)
$$
where
$$
\Upsilon_{i0}(s,t)=n^{-1}\int_s^t\frac{R_{i.}(u) \eta_{1i}(u)
\exp(\bbeta^{\circ T} {\bf
Z}_{i0})}{\{{\mathcal{Y}}(u,\Lambda_0^\circ)\}^{2}}
\sum_{k=1}^{2n}\sum_{l=1}^{m_k}d N_{kl}(u)
$$
and for $j \geq 1$
\begin{eqnarray*}
\hspace*{-0cm}\Upsilon_{ij}(s,t)&=&n^{-1}\int_s^t\frac{ R_{i.}(u)
\eta_{1i}(u)\exp({\bbeta^{\circ T} {\bf
Z}_{ij}})}{\{{\mathcal{Y}}(u,\Lambda_0^\circ)\}^{2}}
\sum_{k=1}^{2n}\sum_{l=1}^{m_k}d N_{kl}(u) \\
& & + \, n^{-1}\sum_{k=1}^{2n}\frac{R_{k.}(s)
\eta_{1k}(s)}{\{{\mathcal{Y}}(s,\Lambda_0^\circ)\}^{2}}
\sum_{l=1}^{m_k}I(T_{kl}>s)\exp({\bbeta^{\circ T} {\bf Z}_{kl}})
\delta_{kl}.
\end{eqnarray*}
The fourth term can be written, by plugging in the representation
for $\tilde\Lambda_0-\Lambda_0^\circ$, as
$$
-n^{-1}\int_0^\tau
\frac{A(s,t)\tilde{p}(s-)}{\tilde{\mathcal{Y}}(s,\Lambda_0^\circ)}
\sum_{i=1}^{2n}\sum_{j=1}^{m_i}I(T_{i0}<s)dM_{ij}(s)
$$
where
$$
A(s,t)=n^{-2}\int_0^t\sum_{k=1}^{2n} \frac{ R_{k.}(s)
\eta_{1k}(s)}
 {\{\mathcal{Y}(s,\Lambda_0^\circ)\}^{2}}\exp(\bbeta^{\circ T} {\bf Z}_{k0})
 \left[\int_s^\tau \{\tilde{p}(v)\}^{-1}d N_{k0}^\star(v)\right]
 \sum_{i=1}^{2n}\sum_{j=1}^{m_i}dN_{ij}(s)
$$
and $N_{k0}^\star(t)=I(T_{k0} \leq t)$. Given all the above, we
get
\begin{eqnarray*}
\lefteqn{\hat{\Lambda}_0(t,\bgamma^\circ)-\Lambda_0^\circ(t)
\approx
 n^{-1}\int_0^t \{{\mathcal{Y}}(s,\Lambda_0^\circ)\}^{-1} \sum_{i=1}^{2n}
 \sum_{j=1}^{m_i} dM_{ij}(s)} \\
&& -n^{-1}\int_0^\tau
\frac{A(s,t)\tilde{p}(s-)}{\tilde{\mathcal{Y}}(s,\Lambda_0^\circ)}
\sum_{i=1}^{2n}\sum_{j=1}^{m_i}I(T_{i0}<s)dM_{ij}(s) \\
&&
-n^{-1}\int_0^t\{\hat{\Lambda}_0(s)-\Lambda_0^\circ(s)\}\sum_{i=1}^{2n}\sum_{j=0}^{m_i}
\Upsilon_{ij}(s,t)d\tilde N_{ij}(s).
\end{eqnarray*}
By solving the above approximation recursively, for the relatives'
failure times, we get
\begin{eqnarray*}
\lefteqn{\hat{\Lambda}_0(t,\bgamma^\circ)-\Lambda_0^\circ(t)
\approx
 \frac{1}{n\hat{p}(t)}\int_0^t\frac{\hat{p}(s-)}{{\mathcal{Y}}(s,\Lambda_0^\circ)} \sum_{i=1}^{2n}
 \sum_{j=1}^{m_i} dM_{ij}(s)} \\
 &&
 +\frac{1}{n\hat{p}(t)}\int_0^\tau B(s,t)\sum_{i=1}^{2n}\sum_{j=1}^{m_i}
 \frac{\tilde{p}(s-)}{{\tilde{\mathcal{Y}}}(s,\Lambda_0^\circ)}I(T_{i0}<s)dM_{ij}(s)\\
&& -\frac{\hat{p}(t-)d N(t)}{n\hat{p}(t)}\int_0^\tau
A(s,t)\sum_{i=1}^{2n}\sum_{j=1}^{m_i}\frac{\tilde{p}(s-)}
{{\tilde{\mathcal{Y}}}(s,\Lambda_0^\circ)}I(T_{i0}<s)dM_{ij}(s)
 \end{eqnarray*}
where $N(t)=\sum_{i=1}^{2n}\sum_{i=1}^{m_i}N_{ij}(t)$,
$$
B(s,t)=n^{-1}\int_0^{t-}A(s,u)\hat{p}(u-)\sum_{i=1}^{2n}\sum_{j=0}^{m_i}\Upsilon_{ij}(u,t-)dN_{ij}(u)
$$
and
$$
\hat{p}(t)=\prod_{s \leq
t}\left[1+\frac{1}{n}\sum_{i=1}^{2n}\sum_{j=0}^{m_i}\Upsilon_{ij}(s,t)d
\tilde{N}_{ij}(s)\right].
$$

\subsection{ Asymptotic normality}
To show that $\hat{\bgamma}$ is asymptotically normally
distributed, we write
\begin{eqnarray}
\lefteqn{{\bf
0}}&=&\bU(\hat{\bgamma},\hat{\Lambda}_0(\cdot,\hat{\bgamma}))\nonumber
\\ &=& \bU(\bgamma^\circ,\Lambda_0^\circ) +
[\bU(\bgamma^\circ,\hat{\Lambda}_0(\cdot,\bgamma^\circ))-\bU(\bgamma^\circ,\Lambda_0^\circ)]
\nonumber\\ & & + \,
[\bU(\hat{\bgamma},\hat{\Lambda}_0(\cdot,\hat{\bgamma}))-
\bU(\bgamma^\circ,\hat{\Lambda}_0(\cdot,\bgamma^\circ))].
\nonumber
\end{eqnarray}
We examine in turn each of the terms on the right-hand side of the above equation.

%\begin{description}

\bsh \underline{Step I} \esh

We can write $\bU(\bgamma^\circ, \Lambda_0^\circ)$ as
$$
\bU(\bgamma^\circ, \Lambda_0^\circ) = \frac{1}{n}
\left(\sum_{i=1}^n
\bxi_i^{(1)}+\sum_{i=1}^{2n}\bxi_i^{(2)}\right).
$$
Here $\bxi_i^{(1)}$ $i=1, \ldots, n$ are iid mean-zero random $(p+1)$-vectors stemming from
the likelihood of the proband data, while $\bxi_i^{(2)}$ $i=1, \ldots, 2n$ are iid mean-zero
random $(p+1)$-vectors stemming from the likelihood of the relatives' data. It follows immediately
from the classical central limit theorem that $n^{-1/2} \bU(\bgamma^\circ, \Lambda_0^\circ)$ is
asymptotically mean-zero multivariate normal.

\bsh \underline{Step II} \esh

Let $\hat{U}_r=U_r(\bgamma^\circ,\hat{\Lambda}_0)$,
$r=1,\ldots,p$, and
$\hat{U}_{p+1}=U_{p+1}(\bgamma^\circ,\hat{\Lambda}_0)$ (in this
segment of the proof, when we write
$(\bgamma^\circ,\hat{\Lambda}_0)$ the intent is to signify
$(\bgamma^\circ,\hat{\Lambda}_0(\cdot,\bgamma^\circ))$. First
order Taylor expansion of $\hat{U}_r$ about $\Lambda_0^\circ$,
$r=1,\ldots,p+1$, gives
$$
\hspace*{-6cm}
n^{1/2}\{U_r(\bgamma^\circ,\hat{\Lambda}_0)-U_r(\bgamma^\circ,{\Lambda}_0^\circ)\}
$$
\begin{equation}\label{eq:tayloru}
\hspace*{1cm} = n^{-1/2} \sum_{i=1}^{2n} \sum_{j=0}^{m_i}
Q_{ijr}(\bgamma^\circ,\Lambda_0^\circ,T_{ij})
\{\hat{\Lambda}_0(T_{ij},\bgamma^\circ)-\Lambda_0^\circ(T_{ij})\}
+ o_p(1),
\end{equation}
where
$$
Q_{ijr}(\bgamma^\circ,\Lambda_0^\circ,T_{ij})=
 \frac{\partial U_r(\bgamma^\circ,\Lambda_0^\circ)}{\partial \Lambda_0^\circ(T_{ij})}
\;\;\;\;\;\; i=1,\ldots,2n \;\; j=0,\ldots,m_i \;\;
r=1,\ldots,p+1.
$$

Now let $N_{ij}^{\star}(t)=I(T_{ij} \leq t)$ $i=1,\ldots,2n$
$j=0,\ldots,m_i$. Then
\begin{equation}
U_r(\bgamma^\circ,\hat{\Lambda}_0)-U_r(\bgamma^\circ,{\Lambda}_0^\circ)
= \frac{1}{n} \sum_{i=1}^{2n} \sum_{j=0}^{m_i}\int_0^{\tau}
Q_{ijr} (\bgamma^\circ,\Lambda^\circ,s)
\{\hat{\Lambda}_0(s,\bgamma^\circ)-\Lambda_0^\circ(s)\}
dN_{ij}^{\star}(s).
\label{zookie}
\end{equation}
Define
$$
{\mathcal A}_r^{(1)}(u)=\frac{\hat{p}(u-)}{\mathcal{Y}(u,\Lambda_0^\circ)}\frac{1}{n}
\sum_{i=1}^{2n}\sum_{j=0}^{m_i}\int_u^\tau
\frac{Q_{ijr}(\bgamma^\circ,\Lambda^\circ,s)}{\hat{p}(s)}dN^{\star}_{ij}(s),
$$
$$
{\mathcal A}_r^{(2)}(u)=\frac{\tilde{p}(u-)}{\tilde{\mathcal{Y}}(u,\Lambda_0^\circ)}\frac{1}{n}
\sum_{i=1}^{2n}\sum_{j=0}^{m_i}\int_0^\tau
\frac{Q_{ijr}(\bgamma^\circ,\Lambda^\circ,s)B(u,s)}{\hat{p}(s)}dN^{\star}_{ij}(s),
$$
and
$$
{\mathcal A}_r^{(3)}(u)=\frac{\tilde{p}(u-)}{\tilde{\mathcal{Y}}(u,\Lambda_0^\circ)}\frac{1}{n}
\sum_{i=1}^{2n}\sum_{j=0}^{m_i}\int_0^\tau
\frac{Q_{ijr}(\bgamma^\circ,\Lambda^\circ,s)\hat{p}(s-)A(u,s)}{\hat{p}(s)}dN_{ij}(s).
$$
Also, let $\alpha_r^{(1)}(u)$, $\alpha_r^{(2)}(u)$ and
$\alpha_{r}^{(3)}(u)$ denote the corresponding limiting values of
${\mathcal A}_r^{(1)}(u)$, ${\mathcal A}_r^{(2)}(u)$ and
${\mathcal A}_r^{(3)}(u)$ as $n$ goes to infinity. Then, after plugging
into (\ref{zookie}) the representation in Sec.~9.3.2 for $\sqrt{n}
[\hat{\Lambda}_0(s,\bgamma^\circ)-\Lambda_0^\circ(s)]$ and replacing the
${\mathcal A}$'s with their limiting values, we obtain
$$
\hspace*{-9cm}
U_r(\bgamma^\circ,\hat{\Lambda}_0)-U_r(\bgamma^\circ,{\Lambda}_0^\circ)
$$
\begin{equation}
\hspace*{1cm} \approx \frac{1}{n} \sum_{k=1}^{2n}
\sum_{l=1}^{m_i}\int_0^{\tau} \left[ \alpha_r^{(1)}(u) +
I(T_{k0}<u)\{\alpha_r^{(2)}(u)-\alpha_{r}^{(3)}(u)\}\right]
 dM_{kl}(u).
\end{equation}
This gives a representation of
$U_r(\bgamma^\circ,\hat{\Lambda}_0)-U_r(\bgamma^\circ,{\Lambda}_0^\circ)$
$r=1,\ldots,p+1$ as the average of independent mean zero iid
random variables. Hence, asymptotic normality follows from the
classical central limit theorem.

\bsh
\underline{Step III}
\esh

First order Taylor expansion of
$\bU(\hat{\bgamma},\hat{\Lambda}_0(\cdot,\hat{\bgamma}))$ about
$\bgamma^\circ=({\bbeta^\circ}^T,\theta^\circ)^T$ gives
 \begin{eqnarray*}
\bU(\hat{\bgamma},\hat{\Lambda}_0(\cdot,\hat{\bgamma}))
=\bU(\bgamma^\circ,\hat{\Lambda}_0(\cdot,\bgamma^\circ))+
\bD(\bgamma^\circ) (\hat{\bgamma}-\bgamma^\circ)^T + o_p(1),
 \end{eqnarray*}
where
 \begin{eqnarray*}
D_{ls}(\bgamma)=\partial
U_l(\bgamma,\hat{\Lambda}_0(\cdot,\bgamma)) / \partial \gamma_s
 \end{eqnarray*}
for $l,s=1,\ldots,p+1$.

Combining the results of Steps I-III above we get that
$n^{1/2}(\hat{\bgamma}-\bgamma^\circ)$ is asymptotically zero-mean
normally distributed with a covariance matrix that can be
consistently estimated by a sandwich-type estimator.

\section{ Acknowledgements}
The authors would like to thank Dr. Kathleen Malone for sharing
the data from the case-control family study of breast cancer,
which motivated the development of this work.  The research was
supported in part by grants from the National Institute of Health
and the United States-Israel Binational Science Foundation (BSF).

\section{ References}
\begin{description}
\item
 Aalen, O. O. (1976). Nonparametric inference in connection with
 multiple decrement models. {\em Scandinavian Journal of Statistics}
 {\bf 3},, 15-27.
\item
 Aalen, O. O. (1992). Modeling heterogeneity in survival analysis
 by the compound Poisson distribution. {\em Annals of Applied Probability}
 {\bf 2,} 951-972.
%\item
%Abramowitz, M. and Stegun, I. A. (Eds.) (1972). {\em
%Handbook of Mathematical Functions with Formulas, Graphs, and
%Mathematical Tables, 9th printing} New York: Dover.
 \item Becher,
H., Schmidt, S., Chang-Claude, J. (2003). Reproductivative factors
and familial predisposition for breast cancer by age 50 years. A
case-control-family study for assessing main effects and possible
gene-environment interaction. {\em International Journal of
Epidemiology} {\bf 32,} 38-48. \item Breslow, N. E. and Day, N. E.
(1980). {\em Statistical methods in cancer research: Vol. 1 - The
analysis of case-control studies.} Lyon, France, IARC Scientific
Publication.
%\item
%Begg, C. B. (2002). On the use of familial aggregation in
%population-based case probands for calculating penetrance. {\em
%Journal of the National Cancer Institute} {\bf 94,} 1221-1226.
\item Clayton, D. G. (1978). A model for association in bivariate
life tables and its application in epidemiological studies of
familial tendency in chronic disease incidence. {\em Biometrika}
{\bf 65,} 141-151.
 \item Coditz, G. A., Rosner, B. A. and Speizer,
F. E. (1996). Risk factors for breast cancer according to family
history of breast cancer. For the Nurses' Health Study Reserch
Group. {\em Journal of National Cancer Institute} {\bf 88,}
365-371.
 \item Fine, J. P., Glidden D. V. and Lee, K. (2003). A
simple estimator for a shared frailty regression model. {\em
Journal of the Royal Statistical Society} {\bf 65}, 317-329.
 \item
Foutz, R. V. (1977). On the unique consistent solution to the
likelihood equation. {Journal of the American Statistical
Association} {\bf 72}, 147-148.
 \item
Genest, C. and MacKay, R. J. (1986). The joy of copulas: Bivariate
distributions with given marginals. {\em The American
Statistician} {\bf 40,} 280-283.
 \item Gill, R. D. (1985).
Discussion of the paper by D. Clayton and J. Cuzick. {\em Journal
of the Royal Statistical Society} {\bf A 148}, 108-109.
 \item
Gill, R. D. (1989). Non- and semi-parametric maximum likelihood
estimators and the Von Mises method (Part 1). {\em Scandinavian
Journal of Statistics}, {\bf 16}, 97-128.
 \item Gill, R. D.
(1992). Marginal partial likelihood. {\em Scandinavian Journal of
Statistics} {\bf 79}, 133-137.
 \item Gorfine, M., Zucker, D. M.
and Hsu, L. (2006). Prospective survival analysis with a general
semiparametric shared frailty model - a pseudo full likelihood
approach {\em Biometrika} {\bf 93}, 735-741.
%\item
%Henderson, R. and Oman, P. (1999). Effect of frailty on marginal
%regression estimates in survival analysis. {\em Journal of the
%Royal Statistical Society} {\bf B 61}, 367-379.
 \item
 Hartman, P. (1973). {\em Ordinary Differential Equations,} 2nd
 ed. (reprinted, 1982), Boston: Birkhauser.
 \item Hopper, J. L., Giels G. G., McCredie, M. R. E., Boyler, P.
(1994). Background, rational and protocol for a case-control
family study of breast cancer. {\em Breast} 79-86.
 \item Hopper,
J. L. (2003). Commentary: Case-control-family design: a paradigm
for future epidemiology research? {\em International Journal of
Epidemiology} {\bf 32,} 48-50.
 \item Hougaard, P. (1986). Survival
models for heterogeneous populations derived from stable
distributions. {\em Biometrika} {\bf 73}, 387-396.
 \item Hougaard,
P. (2000). {\em Analysis of Multivariate Survival data}. New York:
Springer.
 \item Hsu, L., Chen, L., Gorfine, M. and Malone, K.
(2004). Semiparametric estimation of marginal hazard function from
case-control family studies. {\em Biometrics} {\bf 60}, 936-944.
 \item Hsu, L. and Gorfine, M. (2006). Multivariate survival
analysis for case-control family data. {\em Biostatistics} {\bf 7}
387-398.
 \item
 Keiding, N. and Gill, R. (1990). Random truncation models and
 Markov processes. {\em The Annals of Statistics} {\bf 18},
 582-602.
 \item Klein, J. P. (1992). Semiparametric estimation
of random effects using the Cox model based on the EM Algorithm.
{\em Biometrics} {\bf 48}, 795-806.
 \item Kosorok M. R., Lee B. L., Fine J. P. (2004). Robust inference for univariate
proportional hazards regression models. {\em Annals of Statistics}
{\bf 32,} 1448-1491.
 \item Malone, K. E., Daling, J. R., Thompson,
J. D., Cecilia, A. O. Francisco, L. V. and Ostrander E. A. (1998).
BRCA1 mutations and breast cancer in the general population. {\em
Journal of the American Medical Association} {\bf 279,} 922-929.
 \item Malone, K. E., Daling, J. R., Neal, C., Suter, N. M.,
O'brien, C., Cushing-Haugen, K., Jonasdottir, T. J., Thompson, J.
D. and Ostrander E. A. (2000). Frequency of BRCA1/BRCA2 mutations
in a population-based sample of young breast carcinoma cases. {\em
Cancer} {\bf 88,} 1393-1402.
 \item Malone, K. M., Daling, J. R., Doody, D. R., Hsu, L.,
 Bernstein, L., Coates, R. J., Marchbanks, P. A., Simon, M. S.,
 McDonald, J. A., Norman, S. A., Strom, B. L., Burkman, R. T.,
 Ursin, G., Deapen, D., Weiss, L. K., Folger, S., Madeoy, J. J.,
 Friedrichsen, D. M., Suter, N. M., Humphrey, M. C., Spirtas, R.,
 Ostrander, E. A. (2006). Prevalence and predictors of BRCA1 and
 BRCA2 mutations in a population-based study of breast cancer in
 white and black American women aged 35-64 years. {\em Cancer Research}
 {\bf 16}, 8297-8308.
 \item Marshall, A. W. and Olkin, I.
(1988). Families of multivariate distributions. {\em Journal of
the American Statistical Association} {\bf 83,} 834-841.
 \item
McGilchrist, C. A. (1993). REML estimation for survival models
with frailty. {\em Biometrics} {\bf 49}, 221-225.
%\item
%Miki, Y., Swansen, J., Shattuck-Eidens, D., Futreal, P. A.,
%Harshman, K., Tavtigian, S., Lin, Q., Cochran, C., Bennett L. M.,
%Ding, W. (1994). A strong candidate for breast and ovarian cancer
%susceptibility gene BRCA1. {\em Science} {\bf 266,} 66-71.
\item
Nielsen, G. G., Gill, R. D., Andersen, P. K. and Sorensen, T. I.
(1992). A counting process approach to maximum likelihood
estimation of frailty models. {\em Scandinavian Journal of
Statistics} {\bf 19}, 25-43.
\item
Oakes, D. (1989). Bivariate survival models induced by frailties.
{\em Journal of the American Statistical Association} {\bf 84,}
487-493.
\item
Parner, E. (1998). Asymptotic theory for the correlated
gamma-frailty model. {\em Annals of Statistics} {\bf 26}, 183-214.
\item
Prentice, R. L. and Breslow, N. E. (1978). Retrospective studies
and failure time models. {\em Biometrika} {\bf 65,} 153-158.
\item
Ripatti, S. and Palmgren J. (2000). Estimation of multivariate
frailty models using penalized partial likelihood. {\em
Biometrics} {\bf 56}, 1016-1022.
\item
Shih, J. H. and Chatterjee, N. (2002). Analysis of survival data
from case-control family studies. {\em Biometrics} {\bf 58},
502-509.
\item
Shih, J. H. and Louis, T. A. (1995). Inference on the association
parameter in copula models for bivariate survival data. {\em
Biometrics} {\bf 51,} 1384-1399.
%\item
%Struewing, J. P., Hartge, P., Wacholder, S., Backer, S. M.,
%Berlin, M., McAdams, M., Timmerman, M. M., Lawrence, B. C. and
%Tucker, M. A. (1997). The risk of cancer associated with specific
%mutation of BRCA1 and BRCA2 among Ashkenazi Jews. {\em New England
%Journal of Medicine} {\bf 336,} 1401-1408.
\item
Vaida, F. and Xu, R. H. (2000). Proportional hazards model with
random effects. {\em Statistics in Medicine} {\bf 19}, 3309-3324.
%\item
%Wooster, R., Neuhausen, S. L., Mangion, J., Quirk, Y., Ford, D.,
%Collins, N., Nguyen, K., Seal, S., tran, T., Averill, D., Fields,
%P., Marshall, G., Narod, S., Lenoir, G. M., Lynch, H., Feunteun,
%J., Devilee, P., Cornelisse, C. J., Menko, F. H., Daly, P. A.,
%Ormiston, W., McManus, R., Pye, C., Lewis, C. M., CannonAlbright,
%L., Peto, J., Ponder, B. A. J., Skolnick, M. H., Easton, D. F.,
%Goldgar, D. E., Stratton, M. R. (1994). Localization of a breast
%cancer susceptibility gene, BRCA2, to chromosome 13q12-13. {\em
%Science} {\bf 265,} 2088-2090.
%\item Yang, S. and Prentice, R. L. (1999). Semiparametric
%inference in the proportional odds regression model. {\em J. Amer.
%Statist. Ass.} {\bf 94}, 125-136.
 \item Zeger, S., Liang, K.-Y.
and Albert, P. S. (1988). Models for longitudinal data: A
generalized Estimation Equation approach. {\em Biometrics} {\bf
44,} 1049-1060.
 \item Zhao, L. P., Hsu, L., Davidov, O., Potter,
J., Elston, R. C., Prentice, R. L. (1997). Population-based family
study designs. An interdisciplinary research framework for genetic
epidemiology. {\em Genetic Epidemiology} {\bf 14,} 365-388.
 \item Zucker, D. M. (2005). A pseudo partial likelihood method
 for semi-parametric survival regression with covariate errors.
 {\em Journal of the American Statistical Association} {\bf 100}
 1264-1277.
 \item
 Zucker, D. M., Gorfine, M. and Hsu, L. (2006). Pseudo full likelihood
 estimation for prospective survival analysis with a general
 semiparametric shared frailty model: asymptotic theory. {\em
 submitted for publication}.
\end{description}

\newpage
\setlength{\tabcolsep}{2.5mm}
\renewcommand{\baselinestretch}{1.0}
\begin{center}
\begin{table}
\caption{{\em Simulation results: 500 control proband matched with
500 case probands; One relative for each proband; $\beta=0.693$,
$\Lambda_0(t)=t$ , $\theta=2.0$, 500 samples.}}
\begin{center}
\begin{tabular}{lcccccc}
 \hline
 &   \multicolumn{2}{c}{Proposed Method} &
   \multicolumn{2}{c}{Hsu et al.} &
   \multicolumn{2}{c}{Shih and Chatterjjee}    \\
 &   \multicolumn{2}{c}{\rule{4.2cm}{0.05mm}} &
\multicolumn{2}{c}{\rule{4.2cm}{0.05mm}} &
\multicolumn{2}{c}{\rule{4.2cm}{0.05mm}}  \\
   & & Empirical &  & Empirical &  & Empirical \\
  & mean & Standard Error & mean & Standard Error & mean & Standard Error \\
 \hline
 $\hat{\beta}$          &    0.706 & 0.197 & 0.697 & 0.201 & 0.698 & 0.182 \\
 $\hat{\theta}$         &    2.003 & 0.312 & 1.986. & 0.302 & 1.992 & 0.303\\
 $\hat{\Lambda}_0(0.2)$ &    0.201 & 0.034 & 0.204 & 0.030 & 0.202 & 0.029\\
 $\hat{\Lambda}_0(0.4)$ &    0.402 & 0.063 & 0.407 & 0.058 & 0.403 & 0.054\\
 $\hat{\Lambda}_0(0.6)$ &    0.603 & 0.095 & 0.612 & 0.090 & 0.605 & 0.084\\
 $\hat{\Lambda}_0(0.8)$ &    0.809 & 0.136 & 0.820 & 0.131 & 0.811 & 0.122\\
\hline
\end{tabular}
\end{center}
\end{table}
\end{center}

\setlength{\tabcolsep}{2.5mm}
\renewcommand{\baselinestretch}{1.0}
\begin{center}
\begin{table}
\caption{{\em Simulation results for the proposed estimators: 500
control proband matched with 500 case probands; One relative for
each proband; $\Lambda_0(t)=t$; 90\% censoring rate; 2000
samples.}}
\begin{center}
\begin{tabular}{llcccc}
 \hline
  & & & & Empirical & Coverage \\
$\beta$& $\theta$& Estimator &  Mean & Standard Error & Rate  \\
 \hline
0.0   & 2.0 &  $\hat{\beta}$ & -0.013 & 0.217 & 93.5  \\
      &     &  $\hat{\theta}$& 2.127 & 0.872 & 96.0 \\
      &     &  $\hat{\Lambda}_0(0.02)$& 0.020 & 0.006 & 94.2 \\
      &     &  $\hat{\Lambda}_0(0.04)$& 0.041 & 0.010 & 94.7 \\
      &     &  $\hat{\Lambda}_0(0.06)$& 0.061 & 0.015 & 94.8 \\
      &     &  $\hat{\Lambda}_0(0.08)$& 0.081 & 0.020 & 95.0\\
      & 3.0 &  $\hat{\beta}$ & -0.025 & 0.226 & 91.7 \\
      &     &  $\hat{\theta}$& 3.126 & 1.142 & 94.2 \\
      &     &  $\hat{\Lambda}_0(0.02)$& 0.020 & 0.005 & 95.7 \\
      &     &  $\hat{\Lambda}_0(0.04)$& 0.041 & 0.012 & 95.8 \\
      &     &  $\hat{\Lambda}_0(0.06)$& 0.062 & 0.016 & 96.1 \\
      &     &  $\hat{\Lambda}_0(0.08)$& 0.082 & 0.021 & 95.9\\
0.693 & 2.0 &  $\hat{\beta}$ & 0.694 & 0.200 & 96.0 \\
      &     &  $\hat{\theta}$& 2.082 & 0.667 & 94.8 \\
      &     &  $\hat{\Lambda}_0(0.02)$& 0.020 & 0.005 & 95.2 \\
      &     &  $\hat{\Lambda}_0(0.04)$& 0.040 & 0.010 & 95.2 \\
      &     &  $\hat{\Lambda}_0(0.06)$& 0.060 & 0.014 & 96.1 \\
      &     &  $\hat{\Lambda}_0(0.08)$& 0.080 & 0.019 & 96.1\\
      & 3.0 &  $\hat{\beta}$ & 0.689 & 0.206 & 95.4 \\
      &     &  $\hat{\theta}$& 3.172& 0.964 & 95.7 \\
      &     &  $\hat{\Lambda}_0(0.02)$& 0.020 & 0.005 & 94.8 \\
      &     &  $\hat{\Lambda}_0(0.04)$& 0.040 & 0.010 & 95.9 \\
      &     &  $\hat{\Lambda}_0(0.06)$& 0.060 & 0.014 & 96.5 \\
      &     &  $\hat{\Lambda}_0(0.08)$& 0.080 & 0.019 & 95.7\\
 \hline
\end{tabular}
\end{center}
\end{table}
\end{center}

\setlength{\tabcolsep}{6mm}
\renewcommand{\baselinestretch}{1.0}
\begin{center}
\begin{table}
\caption{{\em Simulation results for the proposed estimators: 500
control proband matched with 500 case probands; One relative for
each proband; $\Lambda_0(t)=t$; 60\% censoring rate; 2000
samples.}}
\begin{center}
\begin{tabular}{llcccc}
 \hline
  & & & & Empirical & Coverage \\
$\beta$& $\theta$& Estimator &  Mean & Standard Error & Rate \\
 \hline
0.0   & 2.0 &  $\hat{\beta}$ & 0.007 & 0.191 & 96.0 \\
      &     &  $\hat{\theta}$& 2.031 & 0.348 & 97.5 \\
      &     &  $\hat{\Lambda}_0(0.2)$& 0.200 & 0.035 & 95.1 \\
      &     &  $\hat{\Lambda}_0(0.4)$& 0.399 & 0.067 & 95.1 \\
      &     &  $\hat{\Lambda}_0(0.6)$& 0.598 & 0.099 & 95.1 \\
      &     &  $\hat{\Lambda}_0(0.8)$& 0.797 & 0.135 & 95.1\\
      & 3.0 &  $\hat{\beta}$ & 0.003 & 0.199 & 95.5 \\
      &     &  $\hat{\theta}$& 3.039 & 0.499 & 97.0 \\
      &     &  $\hat{\Lambda}_0(0.2)$& 0.201 & 0.042 & 95.6 \\
      &     &  $\hat{\Lambda}_0(0.4)$& 0.402 & 0.078 & 95.8 \\
      &     &  $\hat{\Lambda}_0(0.6)$& 0.602 & 0.114 & 95.9 \\
      &     &  $\hat{\Lambda}_0(0.8)$& 0.806 & 0.157 & 94.6\\
0.693 & 2.0 &  $\hat{\beta}$ & 0.702 & 0.201 & 96.5 \\
      &     &  $\hat{\theta}$& 2.019 & 0.310 & 96.4 \\
      &     &  $\hat{\Lambda}_0(0.2)$& 0.199 & 0.036 & 95.4 \\
      &     &  $\hat{\Lambda}_0(0.4)$& 0.399 & 0.068 & 96.5 \\
      &     &  $\hat{\Lambda}_0(0.6)$& 0.598 & 0.099 & 96.1 \\
      &     &  $\hat{\Lambda}_0(0.8)$& 0.797 & 0.138 & 95.5\\
      & 3.0 &  $\hat{\beta}$ & 0.699 & 0.211 & 96.5 \\
      &     &  $\hat{\theta}$& 3.037& 0.444 & 97.3 \\
      &     &  $\hat{\Lambda}_0(0.2)$& 0.201 & 0.042 & 95.6 \\
      &     &  $\hat{\Lambda}_0(0.4)$& 0.402 & 0.081 & 96.8 \\
      &     &  $\hat{\Lambda}_0(0.6)$& 0.600 & 0.118 & 95.0 \\
      &     &  $\hat{\Lambda}_0(0.8)$& 0.804 & 0.163 & 93.7\\
 \hline
\end{tabular}
\end{center}
\end{table}
\end{center}

%\newpage
\setlength{\tabcolsep}{6mm}
\renewcommand{\baselinestretch}{1.0}
\begin{center}
\begin{table}
\caption{{\em Simulation results for the proposed estimators: 500
control proband matched with 500 case probands; One relative for
each proband; $\Lambda_0(t)=t$; 30\% censoring rate; 2000
samples.}}
\begin{center}
\begin{tabular}{llcccc}
 \hline
  & & & & Empirical & Coverage \\
$\beta$& $\theta$& Estimator &  Mean & Standard Error & Rate \\
 \hline
0.0   & 2.0 &  $\hat{\beta}$ & 0.007 & 0.047 & 95.5 \\
      &     &  $\hat{\theta}$& 2.013 & 0.247 & 95.3 \\
      &     &  $\hat{\Lambda}_0(0.2)$& 0.200 & 0.037 & 95.5 \\
      &     &  $\hat{\Lambda}_0(0.4)$& 0.397 & 0.073 & 95.0 \\
      &     &  $\hat{\Lambda}_0(0.6)$& 0.596 & 0.110 & 95.1 \\
      &     &  $\hat{\Lambda}_0(0.8)$& 0.794 & 0.147 & 95.5\\
      & 3.0 &  $\hat{\beta}$ & 0.006 & 0.048 & 97.3 \\
      &     &  $\hat{\theta}$& 3.009 & 0.370 & 95.3 \\
      &     &  $\hat{\Lambda}_0(0.2)$& 0.200 & 0.040 & 94.0 \\
      &     &  $\hat{\Lambda}_0(0.4)$& 0.399 & 0.078 & 94.1 \\
      &     &  $\hat{\Lambda}_0(0.6)$& 0.597 & 0.116 & 95.0 \\
      &     &  $\hat{\Lambda}_0(0.8)$& 0.796 & 0.155 & 95.6\\
0.693 & 2.0 &  $\hat{\beta}$ & 0.703 & 0.063 & 96.5 \\
      &     &  $\hat{\theta}$& 1.993 & 0.196 & 95.5 \\
      &     &  $\hat{\Lambda}_0(0.2)$& 0.197 & 0.045 & 94.5 \\
      &     &  $\hat{\Lambda}_0(0.4)$& 0.394 & 0.085 & 94.0 \\
      &     &  $\hat{\Lambda}_0(0.6)$& 0.591 & 0.125 & 94.0 \\
      &     &  $\hat{\Lambda}_0(0.8)$& 0.788 & 0.166 & 94.1\\
      & 3.0 &  $\hat{\beta}$ & 0.703 & 0.061 & 97.2 \\
      &     &  $\hat{\theta}$& 2.999& 0.314 & 96.0 \\
      &     &  $\hat{\Lambda}_0(0.2)$& 0.197 & 0.047 & 94.4 \\
      &     &  $\hat{\Lambda}_0(0.4)$& 0.392 & 0.091 & 94.0 \\
      &     &  $\hat{\Lambda}_0(0.6)$& 0.586 & 0.133 & 94.9 \\
      &     &  $\hat{\Lambda}_0(0.8)$& 0.792 & 0.176 & 95.0\\
 \hline
\end{tabular}
\end{center}
\end{table}
\end{center}

%\newpage
\setlength{\tabcolsep}{2.5mm}
\renewcommand{\baselinestretch}{1.0}
\begin{center}
\begin{table}
\caption{{\em Simulation results of left-restricted data: 500
control proband matched with 500 case probands; One relative for
each proband; $s_0=0.1$, $\beta=0.693$, $\Lambda_0(t)=t$ ,
$\theta=2.0$, 500 samples.}}
\begin{center}
\begin{tabular}{lcccccc}
 \hline
 &   \multicolumn{2}{c}{Proposed Method} &
   \multicolumn{2}{c}{Hsu et al.} &
   \multicolumn{2}{c}{Shih and Chatterjjee}    \\
 &   \multicolumn{2}{c}{\rule{4.2cm}{0.05mm}} &
\multicolumn{2}{c}{\rule{4.2cm}{0.05mm}} &
\multicolumn{2}{c}{\rule{4.2cm}{0.05mm}}  \\
   & & Empirical &  & Empirical &  & Empirical \\
  & mean & Standard Error & mean & Standard Error & mean & Standard Error \\
 \hline
 $\hat{\beta}$          &    0.735 & 0.214 & 0.698 & 0.234 & 0.694 & 0.170 \\
 $\hat{\theta}$         &    2.040 & 0.336 & 2.080. & 0.338 & 2.080 & 0.337\\
 $\hat{\Lambda}_0(0.2)$ &    0.195 & 0.049 & 0.198 & 0.034 & 0.198 & 0.031\\
 $\hat{\Lambda}_0(0.4)$ &    0.392 & 0.090 & 0.402 & 0.068 & 0.401 & 0.062\\
 $\hat{\Lambda}_0(0.6)$ &    0.589 & 0.129 & 0.604 & 0.102 & 0.603 & 0.092\\
 $\hat{\Lambda}_0(0.8)$ &    0.786 & 0.172 & 0.813 & 0.143 & 0.810 & 0.128\\
 $\hat{\Lambda}_0(s_0)$ & 0.098 & 0.025 & - & - & - & - \\
\hline
\end{tabular}
\end{center}
\end{table}
\end{center}

%\newpage
\setlength{\tabcolsep}{2.5mm}
\renewcommand{\baselinestretch}{1.0}
\begin{center}
\begin{table}
\caption{{\em Analysis of a case-control family study of breast
cancer.}}
\begin{center}
\begin{tabular}{lcccccc}
 \hline
 &   \multicolumn{2}{c}{Proposed Method} &
   \multicolumn{2}{c}{Hsu et al.} &
   \multicolumn{2}{c}{Shih and Chatterjee}    \\
 &   \multicolumn{2}{c}{\rule{4.2cm}{0.05mm}} &
\multicolumn{2}{c}{\rule{4.2cm}{0.05mm}} &
\multicolumn{2}{c}{\rule{4.2cm}{0.05mm}}  \\
   & & Bootstrap &  & Bootstrap &  & Bootstrap \\
  & mean & Standard Error & mean & Standard Error & mean & Standard Error \\
 \hline
 $\hat{\beta}$         &   -0.440 & 0.158 & -0.484 & 0.216 & -0.476 & 0.168\\
 $\hat{\theta}$       &     0.952 & 0.443 & 0.889 & 0.443 & 0.944 & 0.460 \\
 $\hat{\Lambda}_0(40)$ &    0.005 & 0.002 & 0.005 & 0.002 & 0.005 & 0.002\\
 $\hat{\Lambda}_0(50)$ &    0.022 & 0.006 & 0.023 & 0.006 & 0.023 & 0.006\\
 $\hat{\Lambda}_0(60)$ &    0.048 & 0.010 & 0.051 & 0.010 & 0.049 & 0.010\\
 $\hat{\Lambda}_0(70)$ &    0.091 & 0.016 & 0.095 & 0.016 & 0.092 & 0.016\\
\hline
\end{tabular}
\end{center}
\end{table}
\end{center}
\end{document}